\documentclass[10pt]{amsart}

\usepackage{graphicx}
\usepackage{epsfig}
\usepackage{graphics}
\usepackage{amsfonts}
\usepackage{amsmath}
\usepackage{latexsym}

\newtheorem{theorem}{Theorem}[section]
\newtheorem{lemma}[theorem]{Lemma}
\newtheorem{proposition}[theorem]{Proposition}
\newtheorem{corollary}[theorem]{Corollary}

\theoremstyle{definition}
\newtheorem{definition}[theorem]{Definition}

\theoremstyle{remark}
\newtheorem{remark}[theorem]{Remark}

\newcommand\smr{S_\mu(\rho,\delta,\alpha,\varepsilon^\rho_M)}
\newcommand\smrhat{\widehat
S_\mu(\rho,\widehat\delta,\alpha,\varepsilon^\rho_M)}

\newcommand\si{\sigma} \newcommand\bproof{\begin{proof}}
\newcommand\eproof{\end{proof}} \newcommand\ra{\rightarrow}
 \newcommand\un{{\relax\ifmmode
1\!\!1\else$1\!\!1$\fi}} \newcommand\ho{H\"{o}lder }
\newcommand\ep{\varepsilon} \newcommand\bs{\backslash}

\newcommand\ijkc{I^c_{j,{\bf k}}}

\newcommand\ijxc{I_{j}^{c}(x)}

\newcommand\kjxc{{\bf k}^c_{j,x}}
\newcommand\kjyc{{\bf k}^c_{j,y}}

\newcommand\supp{\mbox{supp }}
\begin{document}

\title[Heterogeneous ubiquitous systems and Hausdorff
dimension]{Heterogeneous ubiquitous systems in $\mathbb{R}^{d}$ and
Hausdorff dimension}

\author{Julien Barral \and St\'ephane Seuret }

\address{Julien Barral \and St\'ephane Seuret \\ Tel.:
               +33-1-39635279\\ Fax: +33-1-39635995
              \\ INRIA Rocquencourt, BP 105, 78150
               Le Chesnay Cedex, FRANCE}
 \email{julien.barral@inria.fr,
               stephane.seuret@inria.fr} 
\maketitle

\begin{abstract}
Let $\{x_n\}_{n\in\mathbb{N}}$ be a sequence of $[0,1]^d$,
$\{\lambda_n\}_{n\in\mathbb{N}}$ a sequence of positive real numbers
converging to 0, and $\delta>1$. The classical ubiquity results are
concerned with the computation of the Hausdorff dimension of
limsup-sets of the form $S(\delta) = \bigcap_{N\in \mathbb{N}} \bigcup
_{n\geq N} B(x_n,\lambda_n^\delta).$

Let $\mu$ be a positive Borel measure on $[0,1]^d$, $\rho\in (0,1]$
and $\alpha>0$. Consider the finer limsup-set
\[S_{\mu}(\rho,\delta,\alpha)= \bigcap_{N\in \mathbb{N}}\ \ \bigcup
_{n\geq N:\, 
\mu(B(x_n,\lambda^\rho_n)) \sim \lambda_n^{\rho\alpha}}
B(x_n,\lambda_n^\delta).\]

We show that, under suitable assumptions on the measure $\mu$, the
Hausdorff dimension of the sets $S_{\mu}(\rho,\delta,\alpha)$ can be
computed.  Moreover, when $\rho<1$, a yet unknown saturation
phenomenon appears in the computation of the Hausdorff dimension of
$S_{\mu}(\rho,\delta,\alpha)$. Our results apply to several classes of
multifractal measures, and $S(\delta)$ corresponds to the special case
where $\mu$ is a monofractal measure like the Lebesgue measure.

The computation of the dimensions of such sets opens the way to the
study of several new objects and phenomena. Applications are given for
the Diophantine approximation conditioned by (or combined with)
$b$-adic expansion properties, by averages of some Birkhoff sums and
branching random walks, as well as by asymptotic behavior of random
covering numbers.  
\end{abstract}


\section{Introduction}
\label{sec_introduction}

Since the famous result of Jarnik \cite{JARNIK1} concerning
Diophantine approximation and Hausdorff dimension, the following
problem has been widely encountered and studied in various
mathematical situations.

\smallskip

Let $\{x_n\}_{n\in\mathbb{N}}$ be a sequence in a compact metric space
$E$ and $\{\lambda_n\}_{n\in\mathbb{N}}$ a sequence of positive real
numbers converging to 0. Let us define the limsup~set
\[S = \bigcap_{N\in \mathbb{N}} \ \bigcup _{n\geq N}
B(x_n,\lambda_n),\]
and let $D$ be its Hausdorff dimension. Let $\delta>1$. What can be
said about the Hausdorff dimension of the subset $S(\delta)$ of $S$
defined by
\[
S(\delta) = \bigcap_{N\in \mathbb{N}} \ \bigcup _{n\geq N}
B(x_n,\lambda_n^\delta) \,\,\, ?\]

Intuitively one would expect the Hausdorff dimension of $S(\delta)$ to
be lower bounded by ${D}/{\delta}$.  This has been proved to hold
in many cases which can roughly be separated into two classes:
\begin{itemize}
\item
when the sequence $\{(x_n,\lambda_n)\}_n$ forms a sort of ``regular
system'' \cite{BAKERSCHMIDT,BUGEAUD1}, which ensures a strong uniform
repartition of the points $\{x_n\}_n$.
\item
when the sequence $\{(x_n,\lambda_n)\}_n$ forms an ubiquitous system
\cite{Dod,DoMePesVel,JAFF4} with respect to a monofractal measure carried
by the set $S$.
\end{itemize}
Let us mention that similar results are obtained in \cite{STRATURB}
when $E$ is a Julia set.  When $\dim S(\delta) < D$, such subsets
$S(\delta)$ are often referred to as exceptional sets \cite{DODSON1}.
Another type of exceptional sets arises when considering the level
sets of well-chosen functions:
\begin{itemize}
\item
the function associating with each point $x\in [0,1]$ the frequency of
the digit $i\in \{0,1,\ldots,b-1\}$ in the $b$-adic expansion of $x$,
\item
more generally the function associating with each point $x$ the
average of the Birkhoff sums related to some dynamical systems,
\item
the function $x\mapsto h_f(x)$, when $f$ is either a function or a
measure on $\mathbb{R}^d$ and $h_f(x)$ is a measure of the local
regularity (typically an \ho exponent) of $f$ around $x$.
\end{itemize}

It is a natural question to ask whether these two approaches can be
combined to obtain finer exceptional sets.  Let us take an example to
illustrate our purpose.

On one side, it is known since Jarnik's results \cite{JARNIK1} that if
the sequence $\{(x_n,\lambda_n)\}_n $ is made of the rational pairs
$\{( {p}/{q}, {1}/{q^2})\}_{p,q \in \mathbb{N}^{*2},\, p\leq
q}$, then for every $\delta>1$ the subset $S(\delta)$ of $[0,1]$ has a
Hausdorff dimension equal to $ {1}/{\delta}$. In the ubiquity's
setting, this is a consequence of the fact that the family
$\{( {p}/{q}, {1}/{q^2})\}_{p,q\in \mathbb{N}^{*2}}$ forms an
ubiquitous systems associated with the Lebesgue measure
\cite{Dod,DoMePesVel}.

On the other side, given $(\pi_0,\pi_1,\ldots,\pi_{b-1}) \in [0,1]^b$
such that $\sum_{i=0}^{b-1} \pi_i =1$, Besicovitch and later Eggleston
\cite{EGGLESTON} studied the sets $E^{\pi_0,\pi_1,\ldots,\pi_{b-1}}$
of points $x$ such that the frequency of the digit $i\in
\{0,1,\ldots,b-1\}$ in the $b$-adic expansion of $x$ is equal to
$\pi_i$. More precisely, for any $x\in [0,1]$, let us consider the
$b$-adic expansion of $x = \sum_{m=1}^{\infty} x_m b^{-m}$, where
$\forall m$, $x_m \in \{0,1,\ldots,b-1\}$. Let $\phi_{i,n}(x)$ be the
mapping
\begin{equation}
\label{defphi}
 x \mapsto \phi_{i,n}(x) = \frac{\#\{m\leq n: x_m=i\}}{n}.
\end{equation} Then
$E^{\pi_0,\pi_1,\ldots,\pi_{b-1}} = \{x: \forall i\in
\{0,1,\ldots,b-1\}, \, \lim_{n\ra +\infty} \phi_{i,n}(x) =
\pi_i\}$. They found that $\dim E^{\pi_0,\pi_1,\ldots,\pi_{b-1}}=
\sum_{i=0}^{b-1} -\pi_i \log_b \pi_i $.
      
\medskip

We address the problem of the computation of the Hausdorff dimension
of the subsets $E^{\pi_0,\pi_1,\ldots,\pi_{b-1}}_\delta$ of $[0,1]$
defined by
\begin{equation*}
E^{\pi_0,\pi_1,\ldots,\pi_{b-1}}_\delta \!\! = \left\{\! x: \begin{cases}
\begin{array}{c}\exists\,  (p_n,q_n)_n \in (\mathbb{N}^{*2})^{\mathbb{N}}
\mbox { such that }q_n \ra +\infty, \\ \left|x
- {p_n}/{q_n}\right| \leq  {1}/{q_n^{2\delta}} \,\,\mbox{ and }
\forall i \in \{0,\ldots,b-1\}, \\\, \lim_{n \ra +\infty}
\phi_{i,[\log_b (q_n^2)]}\left({p_n}/{q_n}\right) = \pi_i
\end{array}\end{cases} \right\}
\end{equation*}
($[x]$ denotes the integer part of $x$). In other words, we seek in
this example for the Hausdorff dimension of the set of points of
$[0,1]$ which are well-approximated by rational numbers fulfilling a
given Besicovitch condition (i.e.  having given digit frequencies in
their $b$-adic expansion).  This problem is not covered by the works
mentioned above. The main reason is the heterogeneity of the
repartition of the rational numbers satisfying the Besicovitch
conditions.  As a consequence of Theorems \ref{upperbound} and
\ref{theor_lower} of this paper, one obtains
\begin{equation}\label{JaBe}
\dim E^{\pi_0,\pi_1,\ldots,\pi_{b-1}}_\delta = \frac{\sum_{i=0}^{b-1} -\pi_i
\log_b \pi_i }{\delta}.
\end{equation}

\medskip

The key point to achieve this work is to see the Besicovitch condition
as a scaling property derived from a multinomial measure. More
precisely, the computation of the Hausdorff dimensions of the sets
$E^{\pi_0,\pi_1,\ldots,\pi_{b-1}} _\delta$ proves to be a particular
case of the following problem: Let $\mu$ be a positive Borel measure
on the compact metric space $E$ considered above. Given $\alpha>0$ and
$\delta\ge 1$, what is the Hausdorff dimension of the set of points
$x$ of $E$ that are well-approximated by points of
$\{(x_n,\lambda_n)\}_n$ at rate $\delta$, i.e. such that for an
infinite number of integers $n$, $|x-x_n| \leq \lambda_n^\delta$,
conditionally to the fact that the corresponding sequence of couples
$(x_n,\lambda_n)$ satisfies
\begin{equation}
\label{exp}
\displaystyle\lim_{n\to\infty}\frac{\log \mu \big
(B(x_n,\lambda_n)\big )}{\log (\lambda_n)}=\alpha?
\end{equation} 
In other words, if $\varepsilon=(\ep_n)_{n\ge 1}$ is a sequence of
positive numbers converging to 0, what is the Hausdorff dimension of
\begin{equation}
\label{def1}
S_\mu(\delta,\alpha,\varepsilon)=\bigcap_{N\ge 0} \ \ \bigcup_{n\ge N:
\, \lambda_n^{\alpha+\ep_n}\le \mu \left (B(x_n,\lambda_n)\right )\le
\lambda_n^{\alpha-\ep_n}} \hspace{-3mm} B(x_n,\lambda_n^\delta) \ ?
\end{equation}

We study the problem in $\mathbb{R}^d$ ($d\ge 1$). 
An upper bound for the Hausdorff dimension of
$S_\mu(\delta,\alpha,\varepsilon)$ is given by Theorem
\ref{upperbound} for {\em weakly redundant systems} $\{(x_n,
\lambda_n)\}_n$ (see Definition \ref{weaksystem}). Its proof uses
ideas coming from multifractal formalism for measures
\cite{BRMICHPEY,OLSEN}.

Theorem~\ref{theor_lower} (case $\rho=1$) gives a precise lower bound
of the Hausdorff dimension of $S_\mu(\delta,\alpha,\varepsilon)$ when
the family $\{(x_n, \lambda_n)\}_n$ forms a {\em $1$-heterogeneous
ubiquitous system with respect to the measure $\mu$} (see Definition
\ref{deffubi} for this notion, which generalizes the notion of
ubiquitous system mentioned above). It can specifically be applied to
measures $\mu$ that possess some statistical self-similarity property,
and to any family $\{(x_n,\lambda_n)\}_n$ as soon as the support
of~$\mu$ is covered by $\limsup_{n\to\infty}B(x_n,\lambda_n)$.

\medskip

To fix ideas, let us state a corollary of Theorems \ref{upperbound}
and \ref{theor_lower}. This result uses the Legendre transform
$\tau_\mu^*$ of the ``dimension'' function $\tau_\mu$ considered in
the multifractal formalism studied in \cite{BRMICHPEY} (see Section
\ref{subsec_def4} and Definition~\ref{deftaue}).
\begin{theorem}
Let $\mu$ be a multinomial measure on $[0,1]^d$. Suppose that the
fa\-mily $\{(x_n, \lambda_n)\}_n$ forms a weakly redundant
1-heterogeneous ubiquitous system with respect to $\big
(\mu,\alpha,\tau_\mu^*(\alpha)\big )$.

There is a positive sequence $\ep=(\varepsilon_n)_{n\ge 1}$ converging
to 0 at $\infty$ such that
\[\forall\  \delta\ge 1, \ \
\dim\,S_\mu(\delta,\alpha,\varepsilon)= {\tau^*_\mu
(\alpha)}/{\delta}.
\]
\end{theorem}
Examples of remarkable families $\{(x_n,\lambda_n)\}_n$ are discussed
in Section \ref{sec_examples}, as well as examples of suitable
statistically self-similar measures $\mu$.  There, the measures $\mu$
are chosen so that the property (\ref{exp}) has a relevant
interpretation (for instance in terms of the $b$-adic expansion of the
points $x_n$).

\medskip
The formula (\ref{def1}) defining the set
$S_\mu(\delta,\alpha,\varepsilon)$ naturally leads to the question of
conditioned ubiquity into the following more general form: Let $\rho
\in (0,1]$. What is the Hausdorff dimension of
\begin{equation}
\label{def2}
S_\mu(\rho,\delta,\alpha,\varepsilon)=\bigcap_{N\ge 0} \ \
\bigcup_{n\ge N: \, \lambda_n^{\rho(\alpha+\ep_n)}\le \mu \left
(B(x_n,\lambda^\rho_n)\right )\le
\lambda_n^{\rho(\alpha-\ep_n)}}\hspace{-3mm} B(x_n,\lambda_n^\delta) \
?
\end{equation}
Remark that, in (\ref{def1}) and (\ref{def2}), if $\mu$ equals the
Lebesgue measure and if $\alpha=d$, the conditions on
$B(x_n,\lambda_n^\rho)$ are empty, since they are independent of
$x_n$, $\lambda_n$ and $\rho$ (this remains true for a strictly
monofractal measure $\mu$ of index $\alpha$, that is such that
$\exists\, C>0$, $\exists \,r_0$ such that $\forall \,x\in $
supp($\mu$), $\forall\,\, 0<r\leq r_0$, $\, C^{-1} r^\alpha\leq
\mu(B(x,r)) \leq C r^\alpha$).

Again, an upper bound for the Hausdorff dimension of
$S_\mu(\rho,\delta,\alpha,\varepsilon)$ is found in Theorem
\ref{upperbound} for weakly redundant systems.

Theorem~\ref{theor_lower} (case $\rho<1$) yields a lower bound of the
Hausdorff dimension of $S_\mu(\rho,\delta,\alpha,\varepsilon)$ when
$\rho<1$, as soon as the family $\{(x_n,\lambda_n)\}_n$ forms a {\em
$\rho$-heterogeneous ubiquitous system with respect to $\mu$} in the
sense of Definition \ref{deffubir}. The introduction of this dilation
parameter $\rho$ substantially modifies Definition \ref{deffubi} and
the proofs of the results in the initial case $\rho =1$.

\smallskip

As a consequence of Theorem~\ref{theor_lower}, a new
saturation phenomenon occurs for systems that are both weakly
redundant and $\rho$-heterogeneous ubiquitous systems when $\rho
<1$. This points out the heterogeneity introduced when considering
ubiquity conditioned by measures that are not monofractal. The
following result is also a corollary of Theorems \ref{upperbound} and
\ref{theor_lower}.
\begin{theorem}\label{2}
Let $\mu$ be a multinomial measure on $[0,1]^d$. Let $\rho\in
(0,1)$. Suppose that $\{(x_n, \lambda_n)\}_n$ forms a weakly redundant
$\rho$-heterogeneous ubiquitous system with respect to $\big
(\mu,\alpha,\tau_\mu^*(\alpha)\big )$.

There is a positive sequence $\ep=(\varepsilon_n)_{n\ge 1}$ converging
to 0 at $\infty$ such that
\[\forall\ \delta\ge 1, \
 \ \dim\,S_\mu(\rho,\delta,\alpha,\varepsilon)=\min \Big
 (\frac{d(1-\rho)+\rho\tau^*_\mu (\alpha)}{\delta},\tau^*_\mu
 (\alpha)\Big ).
\]
\end{theorem}
Under the assumptions of Theorem~\ref{2}, if $\tau_\mu^*(\alpha)<d$,
although $\delta$ starts to increase from $1$, $\dim\,
S_\mu(\rho,\delta,\alpha,\varepsilon)$ remains constant until $\delta$
reaches the critical value $\frac{d(1-\rho)+\rho\tau^*_\mu
(\alpha)}{\tau^*_\mu (\alpha)}>1$. When $\delta$ becomes larger~than
$\frac{d(1-\rho)+\rho\tau^*_\mu (\alpha)}{\tau^*_\mu (\alpha)}$, the
dimension decreases. This is what we call a saturation~phenomenon.

\medskip

It turns out that conditioned ubiquity as defined in this paper is
closely related to the local regularity properties of some new classes
of functions and measures having dense sets of discontinuities. In
particular, Theorem~\ref{theor_lower} is a determinant tool to analyze
measures constructed as the measures $\nu_{\rho,\gamma,\si}$
\[
\nu_{\rho,\gamma,\si}=\sum_{n\ge 0}\lambda_n^\gamma \, \mu \big
(B(x_n,\lambda^\rho_n)\big )^\sigma\delta_{x_n},
\] where $\delta_{x_n}$ is the probability Dirac mass at $x_n$,
$\rho\in (0,1]$, and $\gamma,\sigma$ are real numbers which make the
series converge. Conditioned ubiquity is also essential to perform the
multifractal analysis of L\'evy processes in multifractal time. These
objects have multifractal properties that were unknown until
now. Their study is achieved in other works \cite{NUGEN,LEVY}.

\medskip

The definitions of weakly redundant and $\rho$-heterogeneous
ubiquitous systems are given in Section \ref{sec_def}. The statements
of the main results (Theorems \ref{upperbound} and \ref{theor_lower})
then~follow. The proofs of Theorem \ref{upperbound}, Theorem
\ref{theor_lower} (case $\rho=1$) and Theorem \ref{theor_lower} (case
$\rho<1$) are respectively achieved in Sections \ref{sec_upper},
\ref{sec_lower} and \ref{sec_lowerrho}. Finally, our results apply to
suitable examples of systems $\{(x_n,\lambda_n)\}_n$ and measures
$\mu$ that are discussed in Section \ref{sec_examples}.

\section{Definitions and statement of results}
\label{sec_def}

It is convenient to endow $\mathbb{R}^d$ with the supremum norm
$\|\cdot\|_{\infty}$ and with the associated distance $
(x,y)\in\mathbb{R}^d\times\mathbb{R}^d\mapsto
\|x-y\|_{\infty}=\max_{1\le i\le d} (|x_i-y_i|)$. All along the paper,
for a set $S$, $|S|$ denotes then the diameter of $S$.
\smallskip

We briefly recall the definition of the generalized Hausdorff measures
and Hausdorff dimension in $\mathbb{R}^d$. Let $\xi$ be a {\em gauge}
function, i.e. a non-negative non-decreasing function on
$\mathbb{R}_+$ such that $\lim _{x\ra 0^+}\xi(x)=0$. Let $S$ be a subset
of $\mathbb{R}^d$. For $\eta>0$, let us  define
\[\mathcal{H}_\eta^\xi(S) = \inf_{\{C_i\}_{i\in \mathcal{I}}:S\subset
\bigcup_{i\in \mathcal{I}} C_i} \ \ \ \sum_{i\in \mathcal{I}}
\xi\left(\left| C_i \right|\right) , \mbox{ (the family $\{C_i\}_{i\in
\mathcal{I}}$ covers $S$)}\] where the infimum is taken over all
countable families $\{C_i\}_{i\in \mathcal{I}} $ such that $\forall
i\in \mathcal{I}$, $|C_i| \leq \eta$. As $\eta$ decreases to 0,
$\mathcal{H}_\eta^\xi(S)$ is non-decreasing, and $\mathcal {H}^\xi(S)
= \lim _{\eta\ra 0} \mathcal{H}_\eta^\xi(S)$ defines a Borel measure
on $\mathbb{R}^d$, called Hausdorff $\xi$-measure.

Defining the family $\xi_\alpha(x)=|x|^\alpha$ ($\alpha \ge 0$), there
exists a unique real number $0\le D\leq d$, called the Hausdorff
dimension of $S$ and denoted $\dim\, S$, such that $\displaystyle D=
\sup\left\{\alpha\ge 0:\mathcal {H}^{\xi_\alpha}(S) = +\infty\right\}
= \inf\left\{\alpha:\mathcal {H}^{\xi_\alpha}(S) =0 \right\}$ (with
the convention $\sup \, \emptyset =0$).
We refer the reader to \cite{MATTILA1,F2} for instance for more
details on Hausdorff dimensions.

\smallskip

Let $\mu$ be a positive Borel measure with a support contained in
$[0,1]^d$. The analysis of the local structure of the measure $\mu$ in
$[0,1]^d$ may be naturally done using a $c$-adic grid ($c\ge 2$). This
is the case for instance for the examples of measures of Section
\ref{sec_examples}. We shall thus need the following definitions.

Let $c$ be an integer $\ge 2$. For every $j\geq 0$,  $\forall {\bf
k}=(k_1,\ldots, k_d) \in \{0,1,\ldots, c^j-1\}^d$, $\ijkc$ denotes the
$c$-adic box $[k_1 c^{-j}, (k_1+1)c^{-j})\times \ldots \times [k_d
c^{-j}, (k_d+1)c^{-j})$.  $\forall x\in [0,1)^d$, $\ijxc$ stands for
the unique $c$-adic box of generation $j$ that contains $x$, and
$\kjxc$ is the unique (multi-)integer such that $\ijxc=I^c_{j,{\bf
k}_{j,x}^c}$. If ${\bf k}=(k_1,\ldots, k_d)$ and ${\bf
k'}=(k'_1,\ldots, k'_d)$ both belong to $\mathbb{N}^d$, $\|{\bf
k}-{\bf k}'\|_{\infty} = \max_i |k_i-k'_i|$. The set of $c$-adic boxes
included in $[0,1)^d$ is denoted by $\mathbf{I}$.

Finally, the lower Hausdorff dimension of $\mu$,
$\underline{\dim}(\mu)$, is classically defined as $\inf\left \{\dim\,
E: E\in\mathcal{B}([0,1]^d),\ \mu(E)>0\right \}$.


\subsection{Weakly redundant systems}
\label{subsec_def2}
Let $\{x_n\}_{n\in\mathbb{N}}$ be a family of points of $[0,1]^d$ and
 $\{\lambda_n\}_{n\in\mathbb{N}}$ a non-increasing sequence of
 positive real numbers converging to 0.  For every $j\ge 0$, let 
\begin{equation}
\label{deftj}
 T_j =\Big \{n:\ 2^{-(j+1)}<\lambda_n\le 2^{-j}\Big \}.
\end{equation}
The following definition introduces a natural property from which an
upper bound for the Hausdorff dimension of limsup-sets (\ref{def1})
and (\ref{def2}) can be derived. {\it Weak redundancy} is slightly
more general than {\it sparsity} of \cite{FALC3}.
\begin{definition}
\label{weaksystem}
The family $\{(x_n,\lambda_n)\}_{n\in\mathbb{N}}$ is said to form a
weakly redundant system if there exists a sequence of integers
$(N_j)_{j\ge 0}$ such that

\noindent
(i) $\lim_{j\to\infty} {\log N_j} / {j}=0$.

\noindent
(ii) for every $j\ge 1$, $T_j$ can be decomposed into $N_j$ pairwise
disjoint subsets (denoted $T_{j,1},\dots,$ $T_{j,N_j} $) such that for
each $1\le i\le N_j $, the family $\big \{B(x_n,\lambda_n):\ n\in
T_{j,i}\big \}$ is composed of disjoint balls.
\end{definition}

One has $\bigcup_{i=1}^{N_j} T_{j,i}=T_j$. Since the $T_{j,i}$ are
pairwise disjoint, any point $x \in [0,1]^d$ is covered by at most
$N_j$ balls $B(x_n,\lambda_n)$, $n\in T_j$. Moreover, for every $i$
and $j$, the number of balls of $T_{j,i}$ is bounded by $C_d 2^{dj}$,
where $C_d$ is a positive constant depending only on $d$. Indeed, if
two integers $n \neq n'$ are such that $\lambda_n$ and $\lambda_{n'}$
belong to $T_{j,i} $, then $\|x_n-x_{n'}\|_\infty \geq 2^{-j}$.

\subsection{Upper bounds for Hausdorff dimensions of conditioned limsup sets}
\label{subsec_def4}

Let $\mu$ be a finite positive Borel measure on $[0,1]^d$. 

We let the reader verify that if $\supp \mu= [0,1]^d$, then the
concave function
\begin{equation}
\label{deftau}
\tau_{\mu,c}:\, q\mapsto
\lim\inf_{j\to\infty}-{j}^{-1}\log_c \sum_{\mathbf{k}\in
\{0,\ldots,c^j-1\}^d}\mu (I^c_{j,\mathbf{k}})^q
\end{equation}
does not depend on the integer $c\ge 2$, and is consequently simply
denoted $\tau_\mu$. This function is considered in the
multifractal formalism for measures of \cite{BRMICHPEY}. Then, the
Legendre transform of $\tau_\mu$ at $\alpha\in \mathbb{R}_+$, denoted
by $\tau_\mu^*$, is defined by
\begin{equation}
\label{deftaue}
\tau_\mu^*:\, \alpha \mapsto\inf_{q\in\mathbb{R}} \big(\alpha q-\tau_\mu
(q)\big)\in\mathbb{R}\cup\{-\infty\}.
\end{equation}

\begin{theorem}\label{upperbound}
Let $\{x_n\}_{n\in\mathbb{N}}$ be a family of points of $[0,1]^d$ and
$\{\lambda_n\}_{n\in\mathbb{N}}$ a non-increasing sequence of positive
real numbers converging to 0. Let $\mu$ be a positive finite Borel
measure with a support equal to $[0,1]^d$. Let $\{\varepsilon_n\}_
{n\in\mathbb{N}}$ be a positive sequence converging to 0, $\rho\in
(0,1]$, $\delta\ge 1$ and $\alpha \ge 0$. Let us define
\[
S_\mu(\rho,\delta,\alpha,\varepsilon)=\bigcap_{N\ge 1} \ \
\bigcup_{n\ge N:\, \lambda_n^{\rho(\alpha+\varepsilon_n)}\le \mu 
(B(x_n,\lambda_n^\rho)  )\le
\lambda_n^{\rho(\alpha-\varepsilon_n)}}
\hspace{-3mm}B(x_n,\lambda_n^\delta).
\]
Suppose that $\{(x_n,\lambda_n)\}_{n \in \mathbb{N}}$ forms a
weakly redundant system. Then
\begin{equation}
\label{def3}
\dim\, S_\mu(\rho,\delta,\alpha,\varepsilon)\le \min \Big
( \ \frac{d(1-\rho)+\rho\tau_\mu^*(\alpha)}{\delta},\tau_\mu^*(\alpha)
\ \Big ).
\end{equation}
Moreover, $S_\mu(\rho,\delta,\alpha,\varepsilon)=\emptyset$ if
$\tau_\mu^*(\alpha)<0$.
\end{theorem}

The result does not depend on the precise value of the
sequence $\{\ep_n\}_n$, as soon as $\lim_{n\ra +\infty} \ep_n=0$. 
The proof of Theorem \ref{upperbound} is given in Section~\ref{sec_upper}.

\subsection{Heterogeneous ubiquitous systems}
\label{subsec_def}

Let $\alpha>0$ and $\beta\in (0,d]$ be two real numbers. They play the
 role respectively of the H\"older exponent of $\mu$ and of the lower
 Hausdorff dimension of an auxiliary measure $m$.

The upper bound obtained by Theorem \ref{upperbound} is rather
natural. Here we seek for conditions that make the inequality
(\ref{def3}) become an equality.  The following Definitions
\ref{deffubi} and \ref{deffubir} provide properties guarantying this
equality.

The notion of {\it heterogeneous ubiquitous system} generalizes the notion
of {\it ubiquitous system} in $\mathbb{R}^d$ considered in \cite{Dod}.
\begin{definition}
\label{deffubi}
The system $\{(x_n,\lambda_n)\}_{n\in\mathbb{N}}$  is
    said to form a $1$-heterogeneous ubiquitous system with respect  to
  $(\mu, \alpha,\beta)$ if conditions {\bf (1-4)} are
  fulfilled.

\smallskip

{\bf (1)} There exist two non-decreasing continuous functions $\phi$
and $\psi$ defined on $\mathbb{R}_+$ with the following properties: \\
- $\varphi (0)=\psi (0)=0$, $r\mapsto r^{-\varphi(r)}$ and $r\mapsto
r^{-\psi(r)}$ are non-increasing near $0^+$, \\ - $\lim_{r\to 0^+}
r^{-\varphi(r)}=+\infty$, and $\forall \, \ep>0$, $r\mapsto
r^{\ep-\varphi (r)}$ is non-decreasing near 0,\\ - $\varphi$ and
$\psi$ verify {\bf (2), (3)} and {\bf (4)}.

\smallskip
{\bf (2)} There exist a measure $m$ with a support equal to $[0,1]^d$ with the
following properties:

\smallskip

\noindent
$\bullet$ $m$-almost every $y\in [0,1]^d$ belongs to $\bigcap_{N \geq
1} \,\,\bigcup_{n\geq N} B(x_n, {\lambda_n}/{2})$, i.e.
\begin{eqnarray} 
&& m \ \Big(\bigcap_{N \geq 1} \,\,\bigcup_{n\geq N} B\big(x_n,
{\lambda_n}/{2} \big) \Big) =\Vert m \Vert.\label {P0}
\end{eqnarray}

\noindent
$\bullet$ One has:
\begin{equation}
\label{P2} 
\begin{cases}
\mbox{\rm{For $m$-almost every $y \in [0,1]^d$, $\exists\ j(y), \
\forall j \geq j(y)$}},\\ \forall\ {\bf k}\mbox{ such that } \|{\bf
k}-\kjyc\|_{\infty} \leq 1, \ {\mathcal P}^1_1(\ijkc) \mbox{ holds},
\end{cases}
\end{equation}
where ${\mathcal P}_{M}^1(I)$ is said to hold for the set $I$ and for the
real number $M\ge 1$~when
\begin{equation}
\label{pm}
{M}^{-1}|I|^{\alpha+\psi(|I|)}\le \mu\big (I \big )\le M
|I|^{\alpha-\psi( |I|)}.
\end{equation}

\smallskip

\noindent
$\bullet$ One has:
\begin{equation}
\label{P1} 
\begin{cases}
\mbox{\rm{For $m$-almost every $y \in [0,1]^d$, $\exists\ j(y), \
\forall j \geq j(y)$}},\\ \forall\ {\bf k}\mbox{ such that } \|{\bf
k}-\kjyc\|_{\infty} \leq 1, \ {\mathcal D}^m_1(\ijkc) \mbox{ holds},
\end{cases}
\end{equation}
where ${\mathcal D}^m_M(I)$ is said to hold for the set $I$ and for
the real number $M>0$~when
\begin{equation}
\label{dm}
m (I )\le M |I|^{\beta-\varphi(|I|)}.
\end{equation}
%

\smallskip

{\bf (3)} (Self-similarity of $m$) For every $c$-adic box $L$ of
    $[0,1)^d$, let $f_L$ denote the canonical affine mapping from $L$
    onto $[0,1)^d$ . There exists a measure $m^L$ on $L$, equivalent
    to the restriction $m_{|L}$ of $m$ to $L$ (in the sense that
    $m_{|L}$ and $m^L$ are absolutely continuous with respect to one
    another), such that property (\ref{P1}) holds for the measure
    $m^L\circ f_L^{-1}$ instead of the measure $m$.

\medskip

For every $n\ge 1$, let us then introduce the sets
\[E_n^L=\left \{x\in L: \ \begin{array}{c} \forall \ j\ge n+\log_c
  \big (|L|^{-1}\big ), \ \forall\ {\bf k} \mbox{ such that } \|{\bf
  k}-\kjxc\|_{\infty} \leq 1,\\ m^L\big (\ijkc\big )\le \left
  (\frac{|\ijkc|}{|L|}\right
  )^{\beta-\varphi\big(\frac{|\ijkc|}{|L|}\big)}\end{array}\right\}.
\] The sets $E_n^L$ form a non-decreasing sequence in $L$, and
by (\ref{P1}) and property {\bf (3)}, $\bigcup_{n\ge 1} E_n^L$ is of
full $m^L$-measure.  One can thus consider the integer
\[
n_L=\inf\left\{n\ge 1:\ m^L(E_n^L)\ge \Vert m^L\Vert/2\right\}.
\]
If $x\in (0,1)^d$ and $j\ge 1$, let us define the set of balls
\[ \mathcal{B}_j(x)=\left \{B(x_n,\lambda_n): x\in B\big
(x_n,{\lambda_n}/{2}\big )\mbox{ and } \lambda_n \in
(c^{-(j+1)},c^{-j}]\right \}.
\] Notice that this set may be empty. Then, if $\delta>1$ and
$B(x_n,\lambda_n)\in \mathcal{B}_j(x)$, consider
$B(x_n,\lambda_n^\delta)$. This ball contains an infinite number of
$c$-adic boxes. Among them, let $\mathbf{B}_n^\delta$ be the set of
$c$-adic boxes of maximal diameter. Then~define 
$$
\displaystyle
B^\delta_j(x)=\bigcup_{B(x_n,\lambda_n)\in
\mathcal{B}_j(x)}\mathbf{B}_n^\delta.
$$
\smallskip
{\bf (4)} (Control of the growth speed $n_L$ and of the mass $\Vert
    m^L\Vert$) There exists a subset $\mathcal{D}$ of $(1,\infty)$
    such that for every $\delta\in \mathcal{D}$, for $m$-almost every
    $x\in \limsup_{n\to\infty}B\big (x_n,{\lambda_n}/{2}\big )$, there
    exists an infinite number of integers $j$ for which there exists
    $L\in \mathcal{B}_j^\delta(x)$ such that
\begin{eqnarray}
n_{L} \le  \log_c \big (|L|^{-1})\varphi (|L|) \ \mbox { and } \
|L|^{\varphi (|L|)}  \le  \Vert m^{L}\Vert.\label{P3}
\end{eqnarray}
\end{definition}

\begin{remark}\label{rem} {\em 1.} {\bf (1)} is a technical assumption.
  In {\bf (2)}, (\ref{P1}) provides a lower bound for the lower
  Hausdorff dimension of the analyzing measure $m$.  (\ref{P2}) yields
  a control of the local behavior of $\mu$, $m$-almost
  everywhere. Then (\ref{P0}) is the natural condition on $m$ to
  analyze ubiquitous properties of $\{(x_n,\lambda_n)\}_n$ conditioned
  by $\mu$.  {\bf (3)} is a kind of self-similarity needed for the
  measure $m$, and {\bf (4)} imposes a control of the growth speed in
  the level sets for the ``copies'' $m^L\circ f_L^{-1}$ of $m$. The
  combination of assumptions {\bf (3)} and {\bf (4)} supplies the
  monofractality property used in classical ubiquity results.

\smallskip

{\em 2.} If $\mu$ is a strictly monofractal measure of exponent $d$
(typically the Lebesgue measure), then {\bf (1-4)} are always
fulfilled with $\alpha=\beta=d$ and $\mu=m$ as soon as (\ref{P0})
holds. In fact, in this case, {\bf (1-4)} imply the conditions required
to be an ubiquitous system in the sense of \cite{Dod,DoMePesVel}.

\smallskip

{\em 3.} For some well-chosen measures $m$, property {\bf (4)}
automatically holds for any system $\{(x_n,\lambda_n)\}_{n\in
\mathbb{N}}$ and $\mathcal{D}=(1,\infty)$. This is due to the fact
that a stronger property holds: {\bf (4')} There exists $J_m$ such
that for every $j\ge J_m$, for every $c$-adic box $L=\ijkc$,
(\ref{P3}) holds.  The first two classes described in
Section~\ref{sec_measures} verify {\bf (4')} (see
\cite{MOIBARRALREFINED}).

The use of the weakened property {\bf (4)} is needed for the last two
examples developed in Section~\ref{sec_measures} and for other
measures constructed similarly (see \cite{MOIBARRALREFINED}). Indeed,
for these kinds of random measures, it was impossible for us to prove
the stronger uniform property {\bf (4')}, and we are only able to
derive (see \cite{MOIBARRALREFINED}) that, with probability 1, {\bf
(4)} holds with a dense countable set $\mathcal{D}$.

\smallskip

{\em 4.} Property {\bf (4)} can be weakened without affecting the
   conclusions of Theorem~\ref{theor_lower} below as follows: {\bf
   weak (4)} There exists a subset $\mathcal{D}$ of $(1,\infty)$ such
   that for every $\delta\in \mathcal{D}$, for $m$-almost every $x\in
   \limsup_{n\to\infty}B\big (x_n,{\lambda_n}/{2}\big )$, there exists
   an increasing sequence $j_k(x)$ such that for every $k$, there
   exists $B(x_{n_k},\lambda_{n_k})\in \mathcal{B}_{j_k(x)}(x)$ as
   well as a $c$-adic box $L_k$ included in
   $B(x_{n_k},\lambda_{n_k}^\delta)$ such that (\ref{P3}) holds with
   $L=L_k$; moreover $\lim_{k\to\infty}\frac{\log |L_k|}{\log
   \lambda_{n_k}}=\delta$.  This weakening, necessary in \cite{LEVY},
   slightly complicates the proof and we decided to only discuss this
   point in this remark.
\end{remark}

In order to treat the case of the limsup-sets (\ref{def2}) defined
with a dilation parameter $\rho<1$, conditions {\bf (2)} and {\bf (4)}
are modified as follows.
\begin{definition}
\label{deffubir}
 Let $\rho <1$. The system $\{(x_n,\lambda_n)\}_{n\in\mathbb{N}}$ is
    said to form a $\rho$-heterogeneous ubiquitous system with respect  to
    $(\mu, \alpha,\beta)$ if the following conditions are
    fulfilled.

\smallskip

{\bf (1)} and {\bf (3)} are the same as in Definition~\ref{deffubi}.

\smallskip

{\bf (2($\rho$))} There exists a measure $m$ with a support equal to
$[0,1]^d$ such~that:

$\bullet$ There exists a non-decreasing continuous function $\chi$
defined on $\mathbb{R}_+$ such that $\chi (0)=0$, $r\mapsto
r^{-\chi(r)}$ is non-increasing near $0^+$, $\lim_{r\to 0^+}r^{-\chi
(r)}=+\infty$, and $\forall \ep,\theta,\gamma>0$, $r\mapsto
r^{\ep-\theta\varphi (r)-\gamma\chi(r)}$ is non-decreasing near 0.

Moreover, for $m$-almost every point $y$, there exists an infinite
number of integers $\{ j_i(y)\}_{i\in\mathbb{N}}$ with the following
property: the ball $B(y, c^{-\rho j_i(y)}) $ contains at least
$c^{j_i(y)(d(1-\rho) -\chi(c^{-j_i(y)}))} $ points $x_n$ such that the
associated couples $( x_n, \lambda_n) $ all satisfy
\begin{eqnarray}\label{P0r} &&\begin{array}{c}
    \lambda_n \in [c^{-j_i(y)+1},
    c^{-j_i(y)(1-\chi(c^{-j_i(y)}))}], \\ \mbox{for every } n'\neq
    n, \ B(x'_n,\lambda'_n)\bigcap B(x_n,\lambda_n) =
    \emptyset.\end{array}
\end{eqnarray}

$\bullet$ (\ref{P2}) and (\ref{P1}) in assumption {\bf (2)} are also
supposed here.

\smallskip

{\bf (4')} There exists $J_m$ such that for every $j\ge J_m$, for
every $c$-adic box $L=\ijkc$, (\ref{P3}) holds. In
particular, {\bf (4)} holds with $ \mathcal{D}=(1,+\infty)$.
\end{definition}

\begin{remark}\label{4} {\em 1.} Heuristically, condition (\ref{P0r})
  ensures that for $m$-almost every $y$, for infinite many numbers
$j$, approximatively $c^{jd(1-\rho)}$ ``disjoint'' couples
$(x_n,\lambda_n)$ such that $\lambda_n \sim c^{-j}$ can be found in
the neighborhood $B(y, c^{-\rho j})$ of $y$. This property is much
stronger than (\ref{P0}).

\smallskip

{\em 2. } Again, the uniform property {\bf (4')} (the same as in item
    {\em 3.} of Remark \ref{rem}) could be weakened into: {\bf
    (4($\rho$))} There exists a subset $\mathcal{D}$ of $(1,\infty)$
    such that for every $\delta\in \mathcal{D}$, for $m$-almost every
    $y$, the sequence $j_i(y)$ of {\bf (2($\rho$))} can be chosen so
    that for every $B(x_n,\lambda_n)$ invoked in (\ref{P0r}), among
    the $c$-adic boxes of maximal diameter $L$ included in
    $B(x_n,\lambda_n^\delta)$, at least one satisfies (\ref{P3}).

Nevertheless, we kept {\bf (4')} because we do not know any example of
system $\{x_n,\lambda_n\}_{n\in\mathbb{N}}$ and of measure $m$ such
that {\bf (2($\rho$))} and the weak form of {\bf (4')} hold but such
that {\bf (2($\rho$))} and {\bf (4')} do not.
\end{remark}

Before stating the results, a last property has to be introduced. Let
$\rho<1$. For every set $I$, for every constant $M>1$,
$\mathcal{P}_{M}^{\rho}(I)$ is said to hold if
\begin{equation}
\label{pmr}
{M}^{-1}|I|^{\alpha +\psi(|I|)) +2 \alpha\chi(|I|))}\le \mu\big (I
\big )\le M |I|^{\alpha-\psi(|I|) -2\alpha\chi(|I|)}.
\end{equation}
 
The dependence in $\rho$ of $\mathcal{P}_{M}^{\rho}(I)$ is hidden in
the function $\chi$ (see (\ref{P0r})).

\smallskip

It is convenient for a $\rho$-heterogeneous ubiquitous system
$\{(x_n,\lambda_n)\}$ ($\rho\in (0,1]$) with respect to
$(\mu,\alpha,\beta)$ to introduce the sequences
$\varepsilon_M^\rho=(\varepsilon_{M,n}^\rho)_{n\ge 1}$ defined for a
constant $M\ge 1$ by $\varepsilon_{M,n}^\rho=\max
(\varepsilon_{M,n}^{\rho,-},\varepsilon_{M,n}^{\rho,+})$, where
\begin{equation}
\label{defep}
 \lambda_n^{\alpha\pm\varepsilon_{M,n}^{\rho,\pm}}=
M^{\mp}(2\lambda_n)^{\alpha\pm\psi(2\lambda_n)\pm
2\alpha\chi(2\lambda_n)} (\mbox{by convention $\chi\equiv 0$ if
$\rho=1$}).
\end{equation}
\subsection{Lower bounds for Hausdorff dimensions of conditioned limsup-sets}
\label{subsec_def3}

The triplets $(\mu,\alpha,\beta)$, together with the auxiliary measure
$m$, have the properties required to study the exceptional sets we
introduced before.

\noindent 
Let $\widehat{\delta}=(\delta_n)_{n\ge
1}\in [1,\infty)^{\mathbb{N}^*}$, $\widetilde\varepsilon=(\varepsilon_n)_{n\ge
1}\in (0,\infty)^{\mathbb{N}^*}$, $\rho\in (0,1]$, $M\ge 1$,~and
\begin{equation}
\label{defgene}
\widehat S_\mu (\rho,\widehat\delta,\alpha,
\widetilde\varepsilon)=\bigcap_{N\ge 1} \ \ \bigcup_{n\ge N:\,
\mathcal{Q}(x_n,\lambda_n,\rho,\alpha,\varepsilon_n)\mbox{ holds}}
B(x_n,\lambda_n^{\delta_n}),
\end{equation}
where $ \mathcal{Q}(x_n,\lambda_n,\rho,\alpha,\varepsilon_n)$ holds
when $\lambda_n^{\rho(\alpha+\varepsilon_n)}\le \mu \big
(B(x_n,\lambda_n^\rho)\big )\le
\lambda_n^{\rho(\alpha-\varepsilon_n)}$. So, if $\widehat{\delta}$ is
a constant sequence equal to some $\delta\ge 1$, the set $\widehat
S_\mu (\rho,\widehat\delta,\alpha,\widetilde\varepsilon) $ coincides
with the set $S_\mu (\rho,\delta,\alpha,\widetilde\varepsilon)$
defined in (\ref{def1}) and considered in~Theorem~\ref{upperbound}.
\smallskip

\begin{theorem}
\label{theor_lower}
Let $\mu$ be a finite positive Borel measure whose support is
  $[0,1]^d$, $\rho \in (0,1]$ and $\alpha,\beta>0$. Let
  $\{x_n\}_{n\in\mathbb{N}}$ be a sequence in $[0,1]^d$ and
  $\{\lambda_n\}_{n\in\mathbb{N}}$ a non-increasing sequence of
  positive real numbers converging to 0.

Suppose that $\{(x_n,\lambda_n)\}_{n\in\mathbb{N}}$ forms a
$\rho$-heterogeneous ubiquitous system with respect to {\bf ($\mu,
\alpha, \beta$)}. Let $\widehat{\mathcal{D}}$ be the set of points
$\delta$ of $\mathbb{R}$ which are limits of a non-decreasing element
of $\big (\{1\}\cup\mathcal{D}\big )^{\mathbb{N}^*}$ (in the case of
$\rho<1$, $\mathcal{D} = (1,+\infty)$).

There exists a constant $M\ge 1$ such that for every $ \delta\in
\widehat{\mathcal{D}}$, one can find a non-decreasing sequence
$\widehat\delta$ converging to $\delta$ and a positive measure
$m_{\rho,\delta}$ which satisfy $m_{\rho,\delta}\left(\smrhat\right)
>0$, and such that for every $x \in \smrhat$, (recall that $\chi\equiv
0$ if $\rho=1$ and the definition of $\varepsilon_M^\rho$
(\ref{defep}))
\begin{eqnarray}
\label{controldemdelta}
&&\limsup_{r\to 0^+} \frac{m_{\rho,\delta }\big(B(x,r)\big
)}{r^{D(\beta,\rho,\delta)-\xi_{\rho,\delta}(r)}}<\infty,\\ \nonumber
\mbox{where } &&\begin{cases}\forall \,\rho\in(0,1],\displaystyle
D(\beta,\rho,\delta) = \min
\Big(\frac{d(1-\rho)+\rho\beta}{\delta},\beta \Big );\\
 \forall \, r>0, \ 
\  \xi_{\rho,\delta}(r) =
(4+d)\varphi(r)+\chi(r).
\end{cases}
\end{eqnarray}

$\widehat\delta$ can be taken equal to the constant sequence
$(\delta)_{n\ge 1}$ if $\delta\in \{1\}\cup \mathcal{D}$.

\end{theorem}

For the two first classes of measures of Section \ref{sec_measures}
(Gibbs measures and products of multinomial measures), {\bf (4')}
holds instead of {\bf (4)} and $\mathcal{D}=(1,+\infty)$, and thus
Theorem \ref{theor_lower} applies with any $\rho\in (0,1]$. To the
contrary, as soon as $\rho<1$, Theorem \ref{theor_lower} does not
apply to the last two classes of Section \ref{sec_measures}
(independent multiplicative cascades and compound Poisson cascades).

\begin{corollary}\label{immediate}
Under the assumptions of Theorem~\ref{theor_lower}, there exists $M\ge
1$ such that for every $ \delta\in \widehat{\mathcal{D}}$, there
exists a non-decreasing sequence $\widehat\delta$ converging to
$\delta$ such that $\mathcal{H}^{\xi_{\rho,\delta}} (\smrhat)>0$.
Moreover, $\widehat\delta=(\delta)_{n\ge 1}$
if $\delta\in \{1\}\cup \mathcal{D}$.

In particular, $\dim\, \smrhat \geq D(\beta,\rho,\delta)$.
\end{corollary}

When $\rho<1$, $\displaystyle D(\beta,\rho,\delta)$ remains constant
and equal to $\beta$ when $\delta $ ranges in $[1,
\frac{d(1-\rho)+\rho \beta}{\beta}]$. This is what we call a
saturation phenomenon. Then, as soon as $\frac{d(1-\rho)+\rho
\beta}{\beta}<\delta$, we are back to a ``normal'' situation where
$\displaystyle D(\beta,\rho,\delta)$ decreases as ${1}/{\delta}$
when $\delta$ increases.

When $\rho=1$, $\displaystyle D(\beta,\rho,\delta)={\beta}/{\delta}$,
thus there is no saturation phenomenon.
\begin{corollary}
\label{exactdim}
Fix $\widetilde \varepsilon=(\varepsilon_n)_{n\ge 1}$ a sequence
 converging to 0 at $\infty$. Under the assumptions of
 Theorem~\ref{upperbound} and Theorem~\ref{theor_lower}, if the family
 $\{(x_n,\lambda_n)\}_{n\in \mathbb{N}}$ both forms a weakly redundant and
 a $\rho$-heterogeneous ubiquitous system with respect to
 $(\mu,\alpha, \tau_\mu^*(\alpha))$, then there exists a constant
 $M\ge 1$ such that for every $\delta \in [\frac{d(1-\rho)+\rho
 \tau_\mu^*(\alpha)}{\tau_\mu^*(\alpha)},+\infty) \cap
 \widehat{\mathcal{D}}$, there exists a non-decreasing sequence
 $\widehat\delta$ converging to $\delta$ such that
\begin{eqnarray*}
\dim \big(\smrhat\big )&=&\dim \big(\smrhat \bs \bigcup_{\delta' >
   \delta} S_\mu(\rho,\delta', \alpha,\widetilde \varepsilon) \big) \\ &=&
   D(\tau_\mu^*(\alpha),\rho,\delta).
\end{eqnarray*}
 Moreover, $\widehat\delta$ can be taken equal to $(\delta)_{n\ge 1}$
   if $\delta\in \{1\}\cup \mathcal{D}$.
\end{corollary}

\begin{remark} {\em 1.} Corollary~\ref{immediate} is an immediate consequence
of Theorem~\ref{theor_lower}.

\smallskip

{\em 2.} In order to prove Corollary~\ref{exactdim}, let us first
 observe that if $\delta>1$ and $\widehat\delta$ is a non-decreasing
 sequence converging to $\delta$ when $n$ tend to $\infty$,
 $\smrhat\subset S_\mu (\rho,\delta',\alpha,\varepsilon^\rho_M)$ for
 all $\delta'<\delta$. Theorem~\ref{upperbound} gives the optimal
 upper bound for $\dim \big(\smrhat\big )$. Again by
 Theorem~\ref{upperbound}, if $\delta\ge \frac{d(1-\rho)+\rho
 \tau_\mu^*(\alpha)}{\tau_\mu^*(\alpha)}$, for $\delta'>\delta$, the
 sets $S_\mu (\rho,\delta',\alpha,\varepsilon^\rho_M)$ form a
 non-increasing family of sets of Hausdorff dimension
 $<D(\tau_\mu^*(\alpha),\rho,\delta)$. This implies
 $\mathcal{H}^{\xi_{\rho,\delta}}\big (\bigcup_{\delta' > \delta}
 S_\mu(\rho,\delta', \alpha,\widetilde \varepsilon) \big)=0$. Then the
 lower bound for $\dim \big(\smrhat \bs \bigcup_{\delta' > \delta}
 S_\mu(\rho,\delta', \alpha,\widetilde \varepsilon) \big)$ is given by
 Corollary~\ref{immediate}. This holds for any sequence $\widetilde
 \varepsilon$ converging to zero.

When $\delta=1$, one necessarily has $\rho=1$ and $\widehat
 \delta=(1)_{n\ge 1}$. The arguments are then similar
 to those used for $\delta>1$.

\smallskip

{\em 3.} The previous statements are still valid if property {\bf
    (4')} is replaced by property {\bf (4($\rho$))} of Remark~\ref{4},
    and in Section~\ref{sec_measures}, the measures considered are
    such that either $\mathcal{D}=(1,\infty)$ or $\mathcal{D}$ is
    dense in $(1,\infty)$.

\end{remark}
\section{Upper bound for the Hausdorff dimension of conditioned limsup-sets: Proof of Theorem~\ref{upperbound}}
\label{sec_upper}

The sequence $\{(x_n,\lambda_n)\}_n$ is fixed, and is supposed to form
a weakly redundant system (Definition \ref{weaksystem}). We shall need
the functions $\forall j\ge 1$
\begin{eqnarray*}
\tau_{\mu,\rho,j}(q) = -j^{-1}\log_2\sum_{n\in T_j}\mu\big
(B(x_n,\lambda_n^\rho)\big )^q \ \mbox{and } \
\tau_{\mu,\rho}(q)  =  \liminf_{j\to\infty} \tau_{\mu,\rho,j}(q),
\end{eqnarray*}
with the convention that the empty sum equals 0 and $\log
(0)=-\infty$.

In the sequel, the Besicovitch's covering theorem  is used repeatedly
\begin{theorem}
\label{besic}(Theorem 2.7 of
\cite{MATTILA1}) Let $d$ be an integer greater than 1. There is a
constant $Q(d)$ depending only on $d$ with the following
properties. Let $A$ be a bounded subset of $\mathbb{R}^d$ and
$\mathcal{F}$ a family of closed balls such that each point of $A$ is
the center of some ball of $\mathcal{F}$.

There are families $\mathcal{F}_1,...,\mathcal{F}_{Q(d)}\subset
\mathcal{F}$ covering $A$ such that each $\mathcal{F}_i$ is disjoint,
i.e.
\[A \subset \bigcup_{i=1}^{Q(d)} \bigcup _{F\in \mathcal{F}_i}F
\mbox{ and }\,\,\forall F,F' \in \mathcal{F}_i \mbox{ with }F\neq
F',\, F\cap F'=\emptyset.\]
\end{theorem}

Let $(N_j)_{j\ge 1}$ be a sequence as in Definition~\ref{weaksystem},
and let us consider for every $j\ge 1$ the associated partition
$\{T_{j,1},\dots,T_{j,N_j}\}$ of $T_j$. For every subset $S$ of $T_j$,
for every $1\le i\le N_j$, Theorem \ref{besic} can be used to extract
from $\big\{B(x_n,\lambda_n^\rho):\ n\in T_{j,i}\cap S\big \}$ $Q(d)$
disjoint families of balls denoted by $T_{j,i,k}(S)$, $1\le k\le
Q(d)$, such that
\begin{equation}
\label{dect}
\bigcup_{n\in T_{j,i}\cap S}B(x_n,\lambda_n^\rho)\subset
\bigcup_{k=1}^{Q(d)} \ \bigcup_{n\in T_{j,i,k}(S)}B(x_n,\lambda_n^\rho).
\end{equation}
Let us then introduce the functions
\begin{eqnarray*}
\widehat\tau_{\mu,\rho,j}(q) & = & -{j}^{-1}\ \log_2 \ \sup_{S\subset
  T_j}\sum_{n\in \bigcup_{i=1}^{N_j} \
  \bigcup_{k=1}^{Q(d)}T_{j,i,k}(S)} \mu\big (B(x_n,\lambda_n^\rho)\big
  )^q\quad (j\ge 1)
\end{eqnarray*}
$\mbox{and } \ \widehat\tau_{\mu,\rho}(q) = \liminf_{j\to\infty}
  \widehat\tau_{\mu,\rho,j}(q)$.  Recall that $\tau_\mu$ is defined in
  (\ref{deftau}).
\begin{lemma}\label{rel}
Under the assumptions of Theorem~\ref{upperbound}, one has
\begin{equation}
\label{def5}
\tau_{\mu,\rho}\ge d(1-\rho)+\rho\tau_\mu\quad
\mbox{and}\quad \widehat\tau_{\mu,\rho}\ge \rho\tau_\mu.
\end{equation}
\end{lemma}

\begin{proof} 
$\bullet$ Let us show the first inequality of (\ref{def5}).

  First suppose that $q\ge 0$. Fix $j\ge 1$ and $1\le i\le N_j$. For
every $n\in T_{j,i}$, $B(x_n,\lambda_n^\rho)\cap[0,1]^d$ is contained
in the union of at most $3^d$ distinct dyadic boxes of generation
$j_\rho:=[j\rho]-1$ denoted $B_1(n),\ldots, B_{3^d}(n)$. Hence
\[\mu\big (B(x_n,\lambda_n^\rho)\big )^q\le \Big(\sum_{i=1}^{3^d}\mu\big
(B_i(n)\big ) \Big )^q \leq 
3^{dq}\sum_{i=1}^{3^d}\mu\big (B_i(n)\big )^q.\]

Moreover, since the balls $B(x_n,\lambda_n)$ ($n\in T_{j,i}$) are
pairwise disjoint and of diameter larger than $2^{-(j+1)}$, there
exists a universal constant $C_d$ depending only on $d$ such that each
dyadic box of generation $j_\rho$ meets less than $C_d2^{d(1-\rho)j}$
of these balls $B(x_n,\lambda_n^\rho)$. Hence when summing over $n\in
T_{j,i}$ the masses $\mu\big (B(x_n,\lambda_n^\rho)\big )^q$, each
dyadic box of generation $j_\rho$ appears at most $C_d2^{d(1-\rho)j}$
times.  This implies that
\begin{eqnarray}
\label{bb}
\sum_{n\in T_{j,i}}\mu\big (B(x_n,\lambda_n^\rho)\big )^q & \le &
3^{dq}C_d 2^{d(1-\rho)j}\sum_{\mathbf{k}\in \{0,\ldots,2^{j_\rho}-1\}^d}\mu
(I_{j,\mathbf{k}})^q\\
\mbox{and }\ \sum_{n\in T_j}\mu\big (B(x_n,\lambda_n^\rho)\big )^q &
\le & 3^{dq}C_d N_j 2^{d(1-\rho)j}
\hspace{-5mm} \sum_{\mathbf{k}\in \{0,\ldots,2^{j_\rho}-1\}^d}\mu
(I_{j,\mathbf{k}})^q.
\end{eqnarray}
Since $\log N_j=o(j)$, one gets $\tau_{\mu,\rho}(q)\ge
d(1-\rho)+\rho\tau_\mu(q)$.

\smallskip

Now suppose that $q<0$. Let us fix $j\ge 1$ and $1\le i\le N_j$. For every
$n\in T_{j,i}$, $B(x_n,\lambda_n^\rho)$ contains a dyadic box $B(n)$
of generation $[j\rho]+1$, and $\mu\big (B(x_n,\lambda_n^\rho)\big
)^q\le \mu\big (B(n)\big )^q$. The same arguments as above also yield
$\tau_{\mu,\rho}(q)\ge d(1-\rho)+\rho\tau_\mu(q)$.

\smallskip

$\bullet$ We now prove the second inequality of (\ref{def5}).

Suppose that $q\ge 0$. Fix $j\ge 1$ and $S$ a subset of $T_j$, as well
as $1\le i\le N_j$ and $1\le k\le Q(d)$. We use the decomposition
(\ref{dect}). Since the balls $B(x_n,\lambda_n^\rho)$ ($n\in
T_{j,i,k}(S)$) are pairwise disjoint and of diameter larger than
$2^{-(j+1)\rho}$, there exists a universal constant $C'_d $, depending
only on $d$, such that each dyadic box of generation $j_\rho$ meets
less than $C'_d $ of these balls. Consequently, the arguments used to
get (\ref{bb}) yield here
\begin{eqnarray*}
\sum_{n\in T_{j,i,k(S)}}\mu\big (B(x_n,\lambda_n^\rho)\big )^q & \le & \,
3^{dq}C'_d \sum_{\mathbf{k}\in
\{0,\ldots,2^{j_\rho}-1\}^d}\hspace{-3mm}\mu (I_{j,\mathbf{k}})^q\\
\mbox{and } \ \sum_{n\in
\bigcup_{i=1}^{N_j} \bigcup_{k=1}^{Q(d)} T_{j,i,k}(S)} 
\hspace{-6mm} \mu\big (B(x_n,\lambda_n^\rho)\big )^q & \le & \,
3^{dq}C'_d Q(d)N_j \sum_{\mathbf{k}\in
\{0,\ldots,2^{j_\rho}-1\}^d} \hspace{-3mm} \mu (I_{j,\mathbf{k}})^q.
\end{eqnarray*}
The right hand side in the previous inequality does not depend on $S$,
hence
\begin{eqnarray*}
\sup_{S\subset T_j}\sum_{n\in
\bigcup_{i=1}^{N_j}\bigcup_{k=1}^{Q(d)}T_{j,i,k}(S)} \hspace{-3mm}
\mu\big (B(x_n,\lambda_n^\rho)\big )^q \le 3^{dq}C'_d Q(d)N_j
\hspace{-2mm} \sum_{\mathbf{k}\in \{0,\ldots,2^{j_\rho}-1\}^d}
\hspace{-3mm} \mu (I_{j,\mathbf{k}})^q,
\end{eqnarray*}
and the conclusion follows.  The case $q<0$ is left to the reader.
\end{proof}

\begin{proof} [Proof of Theorem~\ref{upperbound}]
$\bullet$ {\it First case: $\alpha\le \tau_\mu'(0^-)$}. Hence
$\tau_\mu^*(\alpha)=\inf_{q\ge 0} (\alpha q-\tau_\mu(q))$.  Let us
first prove that $\dim\, S_\mu(\rho,\delta,\alpha)\le
\frac{d(1-\rho)+\rho\tau_\mu^*(\alpha)}{\delta}$.  

Fix $\eta>0$ and $N$ so that $\varepsilon_n<\eta$ for $n\ge N$. Let us
introduce the set $ S_\mu(N,\eta,\rho,\delta,\alpha)= \bigcup_{n\ge N:
\, \lambda_n^{\rho(\alpha+\eta)}\le \mu \big
(B(x_n,\lambda_n^\rho)\big )}B(x_n,\lambda_n^\delta)$. This set is
also~written
\[
S_\mu(N,\eta,\rho,\delta,\alpha)=\bigcup_{j\ge \inf_{n\ge
N}\log_2 (\lambda_n^{-1})} \ \ \bigcup_{n\in
T_j: \, \lambda_n^{\rho(\alpha+\eta)}\le \mu \big
(B(x_n,\lambda_n^\rho)\big
)} \hspace{-3mm}B(x_n,\lambda_n^\delta).
\]
Let us fix $D\ge 0$. Remark that $S_\mu(\rho,\delta,\alpha,\widetilde
\varepsilon) \subset S_\mu(N,\eta,\rho,\delta,\alpha)$. We use this
set as covering of $S_\mu(\rho,\delta,\alpha,\widetilde\varepsilon)$
in order to estimate the $D$-dimensional Hausdorff measure of
$S_\mu(\rho,\delta,\alpha,\widetilde\varepsilon)$.

Fix $q\ge 0$ such that $\tau_\mu(q)>-\infty$. Let $j_q$ be an integer
large enough so that $j\ge j_q$ implies $\tau_{\mu,\rho,j}(q)\ge
\tau_{\mu,\rho}(q)-\eta$. For $j_N=\max\big (j_q,\inf_{n\ge N}\log_2
(\lambda_n^{-1})\big )$, one gets that for some constant $C$ depending
on $D,\delta, \alpha,\eta, \rho$ and $q$,
\begin{eqnarray*}
\mathcal{H}^{\xi_D}_{2\cdot 2^{-j_N\delta}}\big
(S_\mu(\rho,\delta,\alpha,\widetilde\varepsilon)\big)&\le &\sum_{j\ge j_N} \ \
\sum_{n\in T_j: \,\lambda_n^{\rho(\alpha+\eta)}\le \mu \big
(B(x_n,\lambda _n^\rho)\big )}\big |B(x_n,\lambda_n^\delta)\big |^D\\
&\le& \sum_{j\ge j_N} \ \sum_{n\in T_j}
|B(x_n,\lambda_n^\delta)\big
|^D\lambda_n^{-q\rho(\alpha+\eta)}\mu \big (B(x_n,\lambda
_n^\rho)\big )^q\\
&\le & \sum_{j\ge j_N}(22^{-j\delta})^{D}
2^{(j+1)q\rho(\alpha+\eta)}2^{-j\tau_{\mu,\rho,j}(q)}\\
&\le & C\sum_{j\ge j_N} 2^{-j(D\delta -q\rho(\alpha+\eta)
+\tau_{\mu,\rho}(q)-\eta)}.
\end{eqnarray*}
Therefore, if $D> \frac{\rho(\alpha+\eta)
-\tau_{\mu,\rho}(q)+\eta}{\delta}$, $\mathcal{H}^{\xi_D}_{2\cdot
2^{-j_N\delta}}\big (S_\mu(\rho,\delta,\alpha,\widetilde
\varepsilon)\big)$ converges to 0 as $N\to\infty$, and $\dim\,
S_\mu(\rho,\delta,\alpha,\widetilde\varepsilon)\le D$. This yields $\dim\,
S_\mu(\rho,\delta,\alpha,\widetilde \varepsilon)\le \frac{q\rho(\alpha+\eta)
-\tau_{\mu,\rho}(q)+\eta}{\delta}$, which is less than $
\frac{d(1-\rho)+\rho (\alpha q-\tau_\mu(q) )+(q\rho+1)\eta
}{\delta}$ by Lemma~\ref{rel}.  This holds for every $\eta>0$ and
for every $q\ge 0$ such that $\tau_\mu(q)>-\infty$. Finally, $\dim\,
S_\mu(\rho,\delta,\alpha,\widetilde \varepsilon) \le \frac{d(1-\rho)+\rho
\inf_{q\ge 0}\alpha q-\tau_\mu(q)}{\delta} = \frac{d(1-\rho)+\rho
\tau_\mu^*(\alpha)}{\delta}$.

\medskip

Let us now show that $\dim\, S_\mu(\rho,\delta,\alpha, \widetilde
\varepsilon)\le \tau_\mu^*(\alpha)$. This time, for $j\ge 1$ we define
$S_j=\{n\in T_j:\ \lambda_n^{\rho(\alpha+\eta)}\le \mu \big
(B(x_n,\lambda _n^\rho)\big )\}$. By (\ref{dect}), we remark that
\[S_\mu(\rho,\delta,\alpha,\widetilde\varepsilon) \subset \bigcup_{j\ge j_N} \
\bigcup_{i=1}^{N_j} \
\bigcup_{k=1}^{Q(d)} \ \bigcup_{n\in
T_{j,i,k}(S_j)}B(x_n,\lambda_n^\rho).\]
By definition of $\widehat\tau_{\mu,\rho}(q)$, a computation
mimicking the previous one yields
\begin{eqnarray*}
\mathcal{H}^{\xi_D}_{2\cdot 2^{-\rho j_N}}\big
(S_\mu(\rho,\delta,\alpha,\widetilde\varepsilon)\big)&\le &C\sum_{j\ge j_N}
2^{-j(D\rho -q\rho(\alpha+\eta) +\widehat\tau_{\mu,\rho}(q)-\eta)}.
\end{eqnarray*}
Hence $\dim\, S_\mu(\rho,\delta,\alpha,\widetilde\varepsilon) \le
\frac{q\rho(\alpha+\eta) -\widehat\tau_{\mu,\rho}(q)
+\eta}{\rho}$, for every $\eta>0$ and every $q\ge 0$
such that $\tau_\mu(q)>-\infty$. The conclusion follows from
Lemma~\ref{rel}.

Finally, if $\tau_\mu^*(\alpha)<0$ and
$S_\mu(\rho,\delta,\alpha,\widetilde\varepsilon)\neq\emptyset$, the previous
estimates show that $\mathcal{H}^{\xi_D}_{2\cdot 2^{-\rho
j_N}}(S_\mu(\rho,\delta,\alpha,\widetilde\varepsilon))$ is bounded for $D\in
(\tau_\mu^*(\alpha),0)$ (one can formally extend the definition of
$\mathcal{H}^{\xi_D}$ to the case $D<0$). This is a contradiction.

\smallskip

$\bullet$ The proof when $\alpha\ge \tau_\mu'(0^-)$ follows similar
lines.
\end{proof}

\section{Conditioned ubiquity. Proof of Theorem~\ref{theor_lower}
(case $\rho=1$)}
\label{sec_lower}

We assume that a 1-heterogeneous ubiquitous system is fixed. 
With each couple $(x_n, \lambda_n)$ is associated the ball $I_n =
B(x_n, \lambda_n)$. For every $\delta \geq 1$, $I_n^{(\delta)}$
denotes the contracted ball $B(x_n, \lambda_n^\delta)$. The
following property is useful in the sequel. Because of the assumption
{\bf (1)} on $\varphi$ and $\psi$, one has
\begin{equation}
\label{monotone}
\exists C>1, \  \forall \,0<r\le s\le 1,\,\,s^{-\varphi (s)}\le Cr^{-\varphi
(r)} \mbox { and }s^{-\psi (s)}\le Cr^{-\psi (r)}.
\end{equation}
We begin with a simple technical lemma
\begin{lemma}
\label{btol}
Let $y\in [0,1]^d$, and let us assume that there exists an integer
  $j(y)$ such that for some integer $c\geq 2$, (\ref{P2}) and
  (\ref{P1}) hold for $y$ and every $j\geq j(y)$.

There exists a constant $M$ independent of $y$ with the following
  property: for every $n$ such that $y\in B(x_n, {\lambda_n}/{2})$ and
  $\log_c \lambda^{-1}_n \geq j(y)+4$ , ${\mathcal D}^m_{M}(B(y,
  2\lambda_n))$ and ${\mathcal P}^1_{M}(B(x_n, \lambda_n))$ hold.
\end{lemma}
\bproof Let us assume that $y\in B(x_n, {\lambda_n}/{2})$ with
$\lambda_n \leq c^{-j(y)-4}$. Let $j_0$ be the smallest integer $j$
such that $c^{-j}\leq {\lambda_n}/{2}$, and $j_1$ the largest integer
$j$ such that $c^{-j}\geq 2\lambda_n$. One has $j_0 \geq -\log_c
\lambda_n \geq j_1 \geq j(y)$. One thus ensured by construction that $j_0-4
\leq -\log_c \lambda_n \leq j_1 +4$.

Let us recall that $I_{j}(y)$ is the unique $c$-adic box of scale $j$
which contains $y$, and that ${\bf k}_{j,y}$ is the unique $k\in
\mathbb{N}^d$ such that $y \in I_{j,{\bf k}}^c= I_{j}(y)$. One has
$ I_{j_0}^c(y) \subset B(x_n,\lambda_n) \subset
\bigcup_{\|{\bf k} - {\bf k}_{j_1,y}^{c}\|_{\infty} \leq 1}
I_{j_1,{\bf k}}^c$, which yields $\displaystyle \mu(I_{j_0}^c(y))
\leq \mu(B(x_n,\lambda_n)) \leq \sum_{\|{\bf k} - {\bf
k}_{j_1,y}^{c}\|_{\infty} \leq 1} \mu(I_{j_1,{\bf k}}^c)$. Applying
(\ref{P2}) and (\ref{pm}) yields
\[|c^{-j_0}|^{\alpha +\psi(|c^{-j_0}|)} \leq
\mu(B(x_n,\lambda_n)) \leq 3^d |c^{-j_1}|^{\alpha -
\psi(|c^{-j_1}|)}.
\] Combining the fact that $j_0-4 \leq -\log_c \lambda_n \leq j_1 +4$
with (\ref{monotone}) and (\ref{defep}) gives
\[\lambda_n^{\alpha+\ep_{M,n}^{1,+}}= {M}^{-1} |2\lambda_n|^{\alpha +
  \psi(2\lambda_n)} \leq \mu(B(x_n,\lambda_n)) \leq M
|2\lambda_n|^{\alpha
-\psi(2\lambda_n)}=\lambda_n^{\alpha-\ep_{M,n}^{1,-}}
\]
for some constant $M$ that does not depend on $y$.

Similarly, one gets from (\ref{P1}) and (\ref{dm}) that
$\mathcal{D}^m_M\big (B(y,2\lambda_n)\big )$ holds for some constant
$M>0$ that does not depend on $y$.
\eproof

\bproof[of Theorem~\ref{theor_lower} in the case $\rho=1$] All along
the proof, $C$ denotes a constant which depends only on $c$, $\alpha$,
$\beta$, $\delta$, $\varphi$ and $\psi$.

\smallskip

The case $\delta=1$ follows immediately from the assumptions (here
$m_1=m$).

\smallskip

Now let $M\ge 1$ be the constant given by Lemma \ref{btol}.  Let
$\delta \in \widehat{\mathcal{D}}\cap (1,+\infty)$, and let $
\{d_n\}_{n\geq 1}$ be a non-decreasing sequence in $\mathcal{D}$
converging to $\delta$ (if $\delta\in \mathcal{D}$, 
$d_n =\delta$ for every $n$).  For every $k\ge 1$, $j\ge 1$ and $y\in
[0,1]^d$, let
\begin{equation}
\label{njy}
n^{(d_k)}_{j,y}= \inf \left\{n: \lambda_n \leq c^{-j},\ \exists j'
 \geq j: \begin{cases}B(x_n, \lambda_n)\in \mathcal{B}_{j'} (y)\\
 \exists\ L\in\mathbf{B}_n^{d_k},\ \mbox{(\ref{P3}) holds} \end{cases}\!\!\!\!\!\!\!\right\}.
\end{equation}
We shall find a sequence $\widehat \delta=(\delta_j)_{j\ge 1}$,
converging to $\delta$, to construct a generalized Cantor set
$K_\delta$ in $\widehat S_\mu(1,\widehat \delta, \alpha, \ep^1_M)$ and
simultaneously the measure $m_\delta$ on $K_\delta$.  The successive
generations of $c$-adic boxes involved in the construction of
$K_\delta$, namely $G_n$, are obtained by induction.

\medskip

{\bf - First step:} The first generation of boxes defining $K_\delta$
is taken as follows.

Let $L_0=[0,1]^d$. Consider the first element $d_1$ of $\mathcal {D}$
of the sequence converging to $\delta$. We first impose that $\delta_j
:=d_1$, for every $j\geq 1$.

Due to assumptions {\bf (2)}, {\bf (3)} and {\bf (4)}, there exist
$E^{L_0}\subset E^{L_0}_{n_{L_0}}$ such that $m(E^{L_0})\ge \Vert
m\Vert/4$ and an integer $n'_{L_0}\ge n_{L_0}$ such that for all $y\in
E^{L_0}$:

\smallskip
\noindent
- $y \in \bigcap_{N \geq 1} \,\,\bigcup_{n\geq N} B(x_n,
{\lambda_n}/{2})$,

\smallskip
\noindent
- for every
$j\ge n'_{L_0}$, both (\ref{P2}) and (\ref{P1}) hold,

\smallskip
\noindent
- there are infinitely many integers $j$ such that (\ref{P3}) holds
  for some $L\in \mathcal{B}_j^{d_1}(y)$.

\smallskip

In order to construct the first generation of balls of the Cantor set,
we invoke the Besicovitch's covering Theorem \ref{besic}.  We are
going to apply it to $A= E^{L_0}$ and to several families
$\mathcal{F}_1(j)$ of balls constructed as follows.

For $y \in E^{L_0}$, let us denote $n^{(d_1)}_{j,y}$ by $n_{j,y}$.
Then for every $j\geq n'_{L_0}+4$, let us define $\displaystyle
{\mathcal{F}}_1(j)=\left\{B\left(y, 2\lambda_{n_{j,y}}\right):\ y\in
E^{L_0}\right\}$.

\smallskip

The family ${\mathcal{F}}_1(j)$ fulfills the conditions of Theorem
\ref{besic}. Thus, for every $j\geq n'_{L_0}+4$, $Q(d)$ families of
disjoint balls ${\mathcal{F}}^1_1(j),...,{\mathcal{F}}_1^{Q(d)}(j)$,
can be extracted from ${\mathcal{F}}_1(j)$. Therefore, since $m(A)=m(
E^{L_0}) \geq \Vert m\Vert/4$, for some $i$ one has $\displaystyle
m\Big ( \bigcup _{F_{1,k}^i \in {\mathcal{F}}_1^{i}(j)} F_{1,k}^i
\Big) \geq {\Vert m\Vert}/{( 4\,Q(d))}$. Again, one extracts from
${\mathcal{F}}_1^{i}(j)$ a finite family of pairwise disjoint balls $
\widetilde{G}_1(j) = \{ B_1, B_2,\ldots, B_N\}$ such that
\begin{equation}
\label{minormass1}
 m\Big( \bigcup _{ B_k \in \widetilde{G}_1(j)} B_k \Big) \geq
\frac{\Vert m\Vert}{ 8 \,Q(d)}.
\end{equation}

By construction, with each $B_k$ can be associated a point $y_k\in
 E^{L_0} $ so that $ B_k = B(y_k,2\lambda_{n_{j,y_k}})$. Moreover, by
 construction (see (\ref{njy})), $I_{n_{j,y_k}}=
 B(x_{n_{j,y_k}},\lambda_{n_{j,y_k}}) \subset
 B(y_k,2\lambda_{n_{j,y_k}}) = B_k$.  Thus $I_{n_{j,y_k}}^{(d_1)}=
 B(x_{n_{j,y_k}}, \lambda_{n_{j,y_k}}^{d_1})$ is included in
 $B_k$. Finally, Lemma \ref{btol} yield $\mathcal{P}^1_M
 (B(x_{n_{j,y_k}},\lambda_{n_{j,y_k}}))$ and $\mathcal{D}^m_M(B_k)$.

\smallskip

Let $J_k$ be the closure of one of the $c$-adic boxes of maximal
diameter included in $I_{n_{j,y_k}}^{(d_1)}$, and such that both
(\ref{P3}) holds for $J_k$. Such a box exists by
(\ref{njy}). Moreover, by construction one has $|J_k| \leq
|I_{n_{j,y_k}}^{(d_1)}|\leq C|J_k|$ for some universal constant $C$.

 We write $\underline{B_k} = J_k$.  Conversely, if a $c$-adic box $J $
 can be written $\underline{B} $ for some larger ball $B$, one writes
 $B=\overline J$. Therefore, for every closed box $J$ constructed
 above one can ensure by construction that
\begin{equation}
\label{toto1}
{C}^{-1}|J| \leq |\overline J| ^{d_1} \leq C|J|,
\end{equation}
where $C$ depends only on the fixed given sequence $\{d_n\}_n$. We
 eventually~set
\begin{equation}
\label{defg1j}
G_1(j) = \{\underline{B_k}: B_k \in
\widetilde{G}_1(j)\}.
\end{equation}
We notice the following property that will be used in the last step:
By construction, if $I_1$ and $I_2$ belong to $G_1(j)$ then their
distance is at least $\max_{i\in\{1,2\}} ({|\overline{I_i}|}/{2}-
\big ({|\overline{I_i}|}/{2}\big )^{d_1})$, which is larger than
$\max_{i\in\{1,2\}} {|\overline{I_i}|}/{3}$ for $j$ large enough
($d_1>1$ by our assumption).

\smallskip

On the algebra generated by the elements of $G_1(j)$, a probability
measure $m_\delta$ is defined by
\[
m_\delta (I)=\frac{m(\overline I)}{\sum_{J_k\in G_1(j)}m(\overline{J_k})}.
\]

Let $I\in G_1(j)$. By construction, $\mathcal{D}^m_M(\overline{I})$
holds. Using consecutively this fact, (\ref{toto1}) and
(\ref{monotone}), one obtains
\[
 m(\overline I)\le M |\overline I|^{\beta-\varphi (|\overline I|)}\le
C|I|^{{\beta}/{d_1}}|\overline I|^{-\varphi(|\overline
I|)}\le C|I|^{{\beta}/{d_1}}|I|^{-\varphi(| I|)}.
\]
Moreover, by (\ref{minormass1}), and remembering the definition of
$G_1(j)$ (\ref{defg1j}), one gets
\[
\sum_{J_k\in G_1(j)} m(\overline{J_k}) = \sum_{B_k\in
\widetilde{G}_1(j)} m(B_k) \ge \frac{\Vert m\Vert}{8\, Q(d)}.
\]
As a consequence, $\displaystyle \forall \ I\in G_1(j), \,\,m_\delta
(I)\le {8\, Q(d)C}{\Vert m\Vert}^{-1} |I|^{{\beta}/{d_1}}|
I|^{-\varphi(|I|)}$.

By our assumption {\bf (1)}, we can fix $j_1$ large enough so that
\[
\forall \ I\in G_1(j_1), \ {8 \,Q(d)C}{\Vert m\Vert}^{-1} \le |
I|^{-\varphi(|I|)} .
\] We choose the $c$-adic elements of the first generation of the
construction of $K_\delta$ as being those of
$G_1:=G_1(j_1)$. By construction
\begin{equation}
\label{premmaj} \forall \ I\in G_1, \ 
m_\delta (I)\le |I|^{{\beta}/{d_1}-2\varphi(|I|)}.
\end{equation}

One knows that by construction, for every $I \in G_1$, there exists
$y_k \in  {E}^{L_0}$ such that $B(x_{n_{j_1,y_k}},
\lambda_{n_{j_1,y_k}}) \subset \overline I =B(y_k,
2\lambda_{n_{j_1,y_k}})$. 

As a consequence, for every $y\in \bigcup_{I\in G_1}I$, there exists
an integer $n$ such that $\lambda_n \leq c^{-4} $, $|x_n-y|\leq
\lambda_n^{\delta_n}$, and ${\mathcal P}^1_{M}(I_n)={\mathcal
P}^1_{M}(B(x_n,\lambda_n))$ holds.

\medskip
{\bf - Second step:} The second generation of boxes is obtained as
follows. Consider $d_2$, the second element of the sequence
$\{d_n\}_n$ converging to $\delta$. Let $n_1$ be the largest integer
among the $n^{(d_1)}_{j_1,y_k}$, $I\in G_1$. For every $j >n_1$, one
imposes $\delta_j :=d_2$.

\smallskip

Let us focus on one of the $c$-adic boxes $L\in G_1$. The selection
procedure is the same as in the first step.  Due to assumptions {\bf
(2)}, {\bf (3)} and {\bf (4)}, one can find a subset $ E^{L}$ of
$E^{L}_{n^L}$ such that $m^L\big (E^{L}\big )\ge \Vert
m^L\Vert/4$ and an integer $n'_L\ge n_L$ such that for all $y\in E^{L}$:

\smallskip
\noindent - $y\in \bigcap_{N \geq 1} \,\,\bigcup_{n\geq N} B(x_n,
{\lambda_n}/{2})$,

\smallskip
\noindent - $\forall \ j\ge n'_{L}+\log_c \big (|L|^{-1}\big )$,
\begin{equation}
\label{ref4}
 \forall\ {\bf k},\ \|{\bf k}-\kjyc\|_{\infty}\leq 1,\ {\mathcal
D}^{m^L\circ f_L^{-1}}_1\left (f_L(\ijkc)\right)\ \mbox{\rm{and}}\
{\mathcal P}^1_1(\ijkc) \mbox{ hold}.
\end{equation}

\noindent - There are infinitely many integers $j$ such that
(\ref{P3}) holds for some $L\in~\mathcal{B}^{d_2}_j(y)$.

\smallskip

We again apply Theorem \ref{besic} to $A= E^{L}$ and to families
$\mathcal{F}_2(j)$ of balls constructed as above. Hence, for every
$j\geq n'_{L}+\log_c \big (|L|^{-1}\big )+4$, ${\mathcal{F}}_2(j)=
\Big\{B(y,2\lambda_{n^{(d_2)}_{j,y}}):\ y\in E^{L}\Big\}$ ($n^{(d_2)}_{j,y}$ is defined in (\ref{njy})). We set
$n_{j,y}:=n^{(d_2)}_{j,y}$.

The family ${\mathcal{F}}_2(j)$ fulfills the conditions of Theorem
\ref{besic} and covers $ E^{L}$.  By Theorem \ref{besic}, for every
$j\geq n'_{L}+\log_c \big (|L|^{-1}\big )+4$, $Q(d)$ families of
pairwise disjoint boxes $ {\mathcal{F}}^1_2(j), \ldots ,
{\mathcal{F}}_2^{Q(d)}(j)$, whose union covers $ E^{L}$, can be
extracted from $ {\mathcal{F}}_2(j)$. In particular, since
$m^L(A)=m^L( E^{L}) \geq \Vert m^L\Vert/4$, there exists one family of
disjoint boxes $ {\mathcal{F}}_2^{i}(j)=\{F_{2,1}^i, F_{2,2}^i, \ldots
\}$ which satisfies
$ m^L\left( \bigcup _{F_{2,k}^i \in {\mathcal{F}}_2^{i}(j)} F_{2,k}^i
\right) \geq {\Vert m^L\Vert}/{ 4Q(d)}.  $

As in the first step, one  extracts from $\mathcal{F}_2^{i}(j)$ a
finite family of disjoint balls $ \widetilde{G}^L_2(j)= \{ B_1,
B_2,\ldots, B_N\}$ such that
\begin{equation}
\label{minormass2}
m^L\Big( \bigcup _{B_k \in \widetilde{G}^L_2(j)} B_k \Big) \geq
\frac{\Vert m^L \Vert}{ 8 \,Q(d)}.
\end{equation}

As above, with each $B_k$ is associated a point $y_k\in E^{L} $ so
that $B_k = B(y_k , 2\lambda_{n_{j,y_k}})$, and
$I^{(d_2)}_{n_{j,y_k}}\subset I_{n_{j,y_k}}\subset B_k$. Now, notice
that Lemma~\ref{btol} applies with $m^L\circ f^{-1}_L$ instead of $m$
and with the same constant $M$. It follows that $\mathcal{D}^{m^L\circ
f_L^{-1}}_M\big (f_L(B_k)\big )$ and $\mathcal{P}^1_M(I_{n_{j,y_k}})$
hold. Let $J_k$ be the closure of one of the $c$-adic balls of maximal
diameter included in $I_{n_{j,y_k}}^{({d_2})}$ such that (\ref{P3})
holds for $J_k$.

We then define the notation $\underline{B_k} = J_k$, and conversely
 $B_k = \overline{J_k}$.  One also has (\ref{toto1}) (for the same
 constant $C$).  We eventually define
\begin{equation}
\label{defg2j}
G_2^{L}(j) = \{\underline{B_k}: B_k \in
\widetilde{G}^{L}_2(j)\}.
\end{equation}

On the algebra generated by the elements $I$ of $G^L_2(j)$, an
extension of the restriction to the ball $L$ of the measure
$m_\delta$ is defined by
\[
m_\delta (I)=\frac{ m^L(\overline I)}{\sum_{ J_k \in
G^L_2(j)}m^L(\overline {J_k})}\, m_\delta (L).
\]
Let $I\in G_2^{L}(j)$. Since $\mathcal{D}^{m^L\circ f_L^{-1}}_M\big
(f_L(\overline{I})\big )$ holds, one has
\begin{eqnarray*}
m^L(\overline I)&\le & M\left (\frac{|\overline I|}{|L|}\right
)^{\beta-\varphi \big(\frac{|\overline I|}{|L|}\big)}\le
C|I|^{{\beta}/{{d_2}}}|L|^{-\beta}\left (\frac{|\overline
I|}{|L|}\right )^{-\varphi \big(\frac{|\overline I|}{|L|}\big)}\\
&\leq &C|I|^{{\beta}/{{d_2}}}|L|^{-\beta} |I|^{-\varphi (|I|)},
\end{eqnarray*}
where (\ref{monotone}) has been used. Moreover, by (\ref{minormass2})
and (\ref{defg2j}),
\[ \sum_{ J_k \in G^L_2(j)}m^L(\overline {J_k}) = \sum_{ B_k \in
  \widetilde G^L_2(j)}m^L( {B_k})\ge {\Vert
m^L\Vert}/{8 \,Q(d)}.
\]
Consequently, since $m_\delta(L)$ can be bounded using
(\ref{premmaj}), one gets
\begin{eqnarray*}
m_\delta (I) &\le&8 m_\delta(L) {Q(d) }{\Vert m^L\Vert}^{-1}
C|I|^{{\beta}/{{d_2}}}|L|^{-\beta} |I|^{-\varphi (|I|)}\\ & \le&
{8 \,Q(d){\Vert m^L\Vert}^{-1}C|L|^{{\beta}/{{d_1}}-\beta-2\varphi
(|L|) }}|I|^{{\beta}/{{d_2}}-\varphi (|I|)}.
\end{eqnarray*}
By {\bf (1)}, one can choose $j_2(L)$ large enough so that for every
integer $j\ge j_2(L)$, for every $c$-adic ball $I\in G^L_2(j)$,
$\displaystyle {8\, Q(d)C{\Vert m^L\Vert}^{-1}|L|^{{\beta} /
{{d_1}}-\beta-2\varphi (|L|)}} \le |I|^{-\varphi (|I|)}$. Then, taking
$j_2=\max\big \{j_2(L): L\in G_1\big\}$, and defining
\[
G_2=\bigcup_{L\in G_1}G^L_2(j_2),
\] this yields an extension of $m_\delta$ to the algebra generated by
the elements of $G_1\bigcup G_2$ and such that for every $I\in
G_1\bigcup G_2$, $ m_\delta(I)\le |I|^{{\beta}/{{d_2}}-2\varphi
(|I|)}$ (indeed if $I\in G_1$ $|I|^{{\beta}/{d_1}} \leq
|I|^{{\beta}/{d_2}}$ because $d_2\ge d_1$).

Notice that by construction, for every $I\in G_2$, $|I|\le \max_{I\in
G_1} 2(c^{-4}|I|)^{d_2}$.

\medskip

{\bf - Third step:} We end the induction. Assume that $n$ generations
of closed $c$-adic boxes $G_1,\dots, G_n$ are found for some integer
$n\ge 2$. Assume also that a probability measure $m_\delta$ on the
algebra generated by $\bigcup_{1\le p\le n}G_p$ is defined and that
the following properties hold (the fact that this holds for $n=2$
comes from the two previous steps):

\smallskip
{\bf (i)} For every $1\leq p \leq n$, the elements of $G_p$ are closed
pairwise disjoint $c$-adic boxes, and for $2\le p\le n$, $\max_{I\in
G_p}|I|\le 2c^{-4 {d_p}}\max_{I\in G_{p-1}}|I|^{d_p}$. 

For every $1\le p\le n$, with each $I\in G_p$ is associated a ball
$\overline{I}$ such that $I\subset \overline{I}$. There exists a
constant $C>0$ depending on $\{d_n\}_n$ such that $C^{-1}|I|\le
|\overline{I}|^{d_p} \le C|I|$. Moreover, if $I_1$ and $I_2$ belong to
$G_p$ then their distance is at least $\max_{i\in\{1,2\}}
{|\overline{I_i}|}/{2}-~\big ({|\overline{I_i}|}/{2}\big )^{d_p}$.
Moreover, the $\overline{I}$'s ($I\in G_p$) are pairwise disjoint.

\smallskip

{\bf (ii)} For every $2\le p\le n$, each element $I$ of $G_p$ is
included in an element $L$ of $G_{p-1}$. Moreover,
$\overline{I}\subset L$, $\log_c \big (|\overline I|^{-1}\big )\ge
n_{L}+\log_c \big (|L|^{-1}\big )$ and $\overline I\cap
E^{L}_{n_{L}}\neq\emptyset$.

\smallskip

{\bf (iii)} There exists a sequence $\widehat \delta=\{\delta_q\}_{q
  \geq 1}$ such that $\forall \,1\le p\le n$ and $I\in G_p$, there is
  an integer $q$ such that $I\subset I_q^{({\delta_q})} =
  B(x_q,\lambda_q^{\delta_q})\subset \overline{I}$, ${\mathcal
  P}^1_{M}(I_{q})$ holds, and $\delta_q=d_p$. Moreover, the sequence
  $\widehat \delta$ is non-decreasing, and $\forall q$, $\delta_q \leq
  \delta$.

\smallskip

{\bf (iv)} For every $I\in \bigcup_{1\le p  \le n}G_p$, $ m_\delta(I)\le
|I|^{{\beta}/{d_n}-2\varphi (|I|)}$.

\smallskip

{\bf (v)} For every $1\le p\le n-1$, $L\in G_p$, and $I\in G_{p+1}$ such
that $I\subset L$,
\[
m_\delta (I)\le 8\, Q(d) m_\delta(L) \frac{ m^L(\overline I)}{\Vert
m^L\Vert}.
\]


{\bf (vi)} Every $L\in\bigcup_{1\le p \le n}G_p$ satisfies (\ref{P3}).

\smallskip

The constructions of a generation $G_{n+1}$ of $c$-adic balls and an
extension of $m_\delta$ to the algebra generated by the elements of
$\bigcup_{1\le p\le n+1} G_p$ such that properties {\bf (i)} to {\bf
(vi)} hold for $n+1$ are done in the same way as when $n=1$.

\smallskip

By induction, and because of the separation property {\bf (i)}, we
get:\\
\noindent
- a sequence $(G_n)_{n\ge 1}$ and a non-decreasing sequence $\widehat
  \delta$ converging to $\delta$,\\
\noindent
- a probability measure $m_{\delta}$ on $\sigma \big(I: I\in
\bigcup_{n\ge 1}G_n\big )$ \\
\noindent such that properties {\bf (i)} to {\bf (vi)} hold for every
$n\ge 2$. We now define
\[
K_\delta=\bigcap_{n\ge 1} \, \bigcup_{I\in G_n}I.
\]

By construction, $m_\delta(K_\delta)=1$ and because of property {\bf
(iii)}, one has $K_\delta\subset\widehat S_\mu(1,\widehat \delta,
\alpha, \ep^1_M)$.  The measure $m_{\delta}$ can be extended to
$\mathcal{B}([0,1]^d)$ by the usual way: $m_\delta(B):=m_\delta(B\cap
K_\delta)$ for $B\in \mathcal{B}([0,1]^d)$. Finally, since
$\delta_n\leq \delta$ for every $n\geq 1$, property {\bf (iv)} implies
that for every $I\in \bigcup_{p \geq 1}G_p$,
\begin{equation}
\label{new4}
m_\delta(I)\le |I|^{{\beta}/{\delta}-2\varphi (|I|)}.
\end{equation}

{\bf - Last step:} Proof of (\ref{controldemdelta}). If $I \in G_n$,
we set $g(I) =n $.

Let us fix $B$ an open ball of $[0,1]^d$ of length less than the one
of the elements of $G_1$, and assume that $B\cap K_\delta \neq
\emptyset$. Let $L$ be the element of largest diameter in
$\bigcup_{n\ge 1}G_n$ such that $B$ intersects at least two elements
of $G_{g(L)+1}$ included in $L$. Remark that this implies that $B$
does not intersect any other element of $G_{g(L)}$, and as a
consequence $ m_\delta (B)\le m_\delta (L)$.

\smallskip

Let us distinguish three cases:

\smallskip
\noindent
$\bullet$ If $|B|\ge |L|$, one has by (\ref{new4})
\begin{equation}
\label{majj1}
m_\delta (B)\le m_\delta (L)\le |L|^{{\beta}/{\delta}-2\varphi
(|L|)} \le C |B|^{{\beta}/{\delta}-2\varphi (|B|)}.
\end{equation}

\smallskip
\noindent
$\bullet$ If $|B|\leq c^{-n_{L}-3}|L|$, let $L_1,\dots,L_p$ be the
elements of $G_{g(L)+1}$ that intersect $B$. We use property {\bf (v)}
to get
\begin{equation}
\label{toto2}
m_\delta(B) = \sum_{i=1}^p m_\delta (B\cap L_i) \le m_\delta(L)
\frac{8 \,Q(d)}{\Vert m^L\Vert}\sum_{i=1}^p m^L(\overline L_i).
\end{equation}
Let $j_0$ be the unique integer such that $c^{-j_0}\le |B| <
c^{-j_0+1}$.  Assume $B$ intersects for instance the boxes $L_{i_1}$
and $L_{i_2}$. Then, by {\bf (i)}, one has $|B| \geq {\max(|\overline
L_{i_1}|,|\overline L_{i_2}|) }/{3}$ when $j_0$ is large enough.
Hence, if $|B|$ is small enough, one has $|B| \geq {(\max_{i=1,...,p}
|\overline L_i|)}/{3}$ and the scale of the boxes $\overline L_i$
(defined as $[-\log_c |\overline L_i|] $) is always larger than
$j_0-[\log_c 3]\geq j_0 -2$.

\medskip

By property {\bf (ii)}, for each $i\in\{1,\dots d\}$, one has
$E^L_{n_L}\cap \overline L_i \neq \emptyset$. Let $y \in E^L_{n_L}\cap
\overline L_i$ for some $i$, and let us consider the $c$-adic box
$I_{j_0-2,{\bf k}_{j_0-2,y}}^c$. For every $z\in \overline L_i$,
$|y-z| \leq c^{-(j_0-2)}$. One deduces that
\[\overline L_i \subset \bigcup_{{\bf k}:\ \|{\bf k}-{\bf
k}_{j_0-2,y}\|_{\infty} \leq 1} I_{j_0-2,{\bf k}}^c.\]

The ball $B$ intersects $L_i$, thus the distance between $y$ and $B$
is at most $c^{-(j_0-2)}$.  As a consequence, if $L_{i'}\neq L_i$, the
distance between $y$ and $L_{i'}$ is lower than $c^{-(j_0-3)}$. This
implies that
\begin{equation}
\label{union}
\bigcup_{i=1}^p \overline L_i \subset \bigcup_{{\bf k}:\ \|{\bf
k}-{\bf k}_{j_0-3,y}\| _{\infty}\leq 1} I_{j_0-3,{\bf k}}^c.
\end{equation}

Since $y\in E^L_{n_L}$ and $j_0\geq -\log_c |L| +n_L+3$, assumption {\bf
  (3)} ensures the control of the $m$-mass of the unions of all the
balls that appear on the left hand-side of (\ref{union}) by the sum of
the masses of the $3^d$ $c$-adic boxes $I^c_{j_0-3,{\bf k}}$, $\|{\bf
  k}-{\bf k}_{j_0-3,y}\| _{\infty}\leq 1$. These boxes all satisfy
\[ m^L(I^c_{j_0-3,{\bf k}}) \leq \left(\frac{|I^c_{j_0-3,{\bf
k}}|}{|L|}\right )^{\beta-\varphi \big(\frac{|I_{j_0-3,{\bf
k}}^c|}{|L|}\big)}\le C\left (\frac{|B|}{|L|}\right
)^{\beta}\left(\frac{|B|}{|L|}\right)^{-\varphi\big(\frac{|B|}{|L|}\big)}
\]
where $C$ depends only on $\beta$. Injecting this in (\ref{toto2}) and
using that the $\overline L_i$ are pairwise disjoint, one obtains that for
$|B|$ small enough 
\begin{eqnarray*}
m_\delta (B)&\le &m_\delta(L) \frac{8 \, Q(d)}{\Vert
m^L\Vert}\sum_{i=1}^p m^L(\overline L_i)\\ &\leq& m_\delta (L) \frac{8
\,Q(d)}{\Vert m^L\Vert}3^d\, C \left (\frac{|B|}{|L|}\right
)^{\beta}\left (\frac{|B|}{|L|}\right
)^{-\varphi\big(\frac{|B|}{|L|}\big)}\\
&\le &m_\delta (L) \frac{C}{\Vert
m^L\Vert} \left (\frac{|B|}{|L|}\right )^{\beta} |B|^{-\varphi (B)},
\end{eqnarray*}
where $C$ takes into account all the constant factors. We then use
consecutively two facts.  First, by (\ref{new4}), $\displaystyle
m_\delta (L)\le |L|^{{\beta}/{\delta}}|L|^{-2\varphi(|L|)}\le C
|L|^{{\beta}/{\delta}}|B|^{-2\varphi(|B|)}$, which implies, since
$r\mapsto r^{\beta(1-1/\delta)}$ is bounded near 0,
\[
m_\delta (B)\le \frac{C}{\Vert
m^L\Vert}|B|^{{\beta}/{\delta}}|B|^{-3\varphi(|B|)}\left
(\frac{|B|}{|L|}\right )^{\beta(1-1/\delta)}\le \frac{C}{\Vert
m^L\Vert}|B|^{{\beta}/{\delta}}|B|^{-3\varphi(|B|)}.
\]

Second, {\bf (vi)} allows to upper bound ${\Vert m^L\Vert}^{-1}$ by
$|L|^{-\varphi(L)}$, which yields
\begin{equation}\label{cas1}
m_\delta (B)\le
C|L|^{-\varphi(|L|)}|B|^{{\beta}/{\delta}}|B|^{-3\varphi(|B|)}\le
C|B|^{{\beta}/{\delta}}|B|^{-4\varphi(|B|)}.
\end{equation}

\smallskip

$\bullet$ $c^{-n_L-3}|L|< |B|\le |L|$: one needs at most
$c^{d(n_L+4)}$ contiguous boxes of diameter $c^{-n_L-3}|L|$ to cover
$B$. For these boxes, the estimate (\ref{cas1}) can be used. Also one
knows by {\bf (vi)} that $c^{n_L} \leq |L|^{-\varphi(L)}$, so for
$|B|$ small enough
\begin{eqnarray*}
&& m_\delta (B)\le C c^{d(n_L+4)}\big (c^{-n_L-3}|L|\big
)^{{\beta}/{\delta}-4\varphi(c^{-n_L-3}|L|)}\le C c^{d
n_L}|B|^{{\beta}/{\delta}-4\varphi(|B|)} \\ &&\le C|L|^{-d\varphi
(|L|)}|B|^{{\beta}/{\delta}-4\varphi(|B|)}\le C
|B|^{{\beta}/{\delta}-(4+d)\varphi(|B|)}.
\end{eqnarray*}
Remembering (\ref{majj1}) and (\ref{cas1}), and using assumption {\bf
(1)}, one gets a constant $C$ such that for every non-trivial ball $B$
of $[0,1]^d$ small enough, one has $m_\delta (B)\le C
|B|^{{\beta}/{\delta}} |B|^{- (4+d)\varphi(|B|)}$. This yields
(\ref{controldemdelta}).
\eproof
\section{Dilation and Saturation. Proof of Theorem~\ref{theor_lower}
(Case $\rho<1$)}
\label{sec_lowerrho}

The introduction of the condition (\ref{P0r}) induces a modification
in the construction of the Cantor set with respect to the case
$\rho=1$, in the selection of the couples $(x_n,\lambda_n)$.  The
following lemma is comparable with Lemma \ref{btol}
\begin{lemma}
\label{btolr}
Let $y\in [0,1]^d$, and assume that (\ref{P2}) and (\ref{P1}) hold for
  $y$ when $j\geq j(y)$ for some integer $j(y)$.  There exists a
  constant $M$ independent of $y$ with the following property: for
  every integer $j$ such that $j(1-\chi (c^{-j}) )\geq
  \frac{j(y)+5}{\rho}$, for every integer $n$ such that $\lambda_n \in
  [c^{-j+1}, c^{-j(1-\chi(c^{-j}))}]$ and
\begin{equation}
\label{wizr}
B(y, (c^{\rho}-1)c^{-j\rho}) \subset B(x_n, \lambda_n^\rho) \subset
   B(y,c^{-j\rho(1-\chi(c^{-j}))}),
\end{equation}
then ${\mathcal P}^\rho_{M}(B(x_n, \lambda_n^\rho))$ holds. Moreover,
the same constant $M$ can be chosen so that ${\mathcal D}^m_{M}(B(y,
r))$ holds for $r\in (0, c^{-j(y)-1})$.
\end{lemma}

\bproof Let us fix $j$ such that (\ref{wizr}) holds, and let us denote
$j_1$ the integer $[j\rho]+2$ and $j_2$ the integer
$[j\rho(1-\chi(c^{-j}))]-2$.  By definition of $j_1$ and $j_2$,
(\ref{wizr}) implies that $I_{j_1}^c(y)\subset B(x_n,
\lambda_n^\rho) \subset \bigcup_{\|{\bf k} - {\bf
k}_{j_2,y}^c\|_{\infty} \leq 1} I_{j_2,{\bf k}}^c$. Combining this
with (\ref{P2}) yields
\begin{equation}
\label{wizzz}
(c^{-j_1})^{\alpha+\psi(c^{-j_1})} \leq \mu(B(x_n,
\lambda_n^\rho)) \leq 3^d (c^{-j_2})^{\alpha-\psi(c^{-j_2})}.
\end{equation}
 One has $c^{-j_1}\leq 2\lambda_n^\rho =|B(x_n,
\lambda_n^\rho)|\leq 2c^{-j_2}$, but by (\ref{wizr}) one also has
\begin{equation}
\label{www}
{C}^{-1} (2c^{-j_2})^{\frac{1}{1-\chi(c^{-j})}}\leq
2\lambda_n^\rho\leq C(2c^{-j_1})^{1-\chi(c^{-j})}
\end{equation}
for some constant $C$ independent of $y$ and $j$.  Hence, using the
monotonicity of $r\mapsto r^{-\psi(r)}$, (\ref{wizzz}) and (\ref{www})
yields the two inequalities
\begin{eqnarray*}
M^{-1}(2\lambda_n^\rho)^{\frac{\alpha}{1-\chi(c^{-j})}}
(2\lambda_n^{\frac{\rho}{1-\chi(c^{-j})}})^
{\psi\big(2\lambda_n^{\frac{\rho}{1-\chi(c^{-j})}}\big)} & \leq &
\mu(B(x_n,\lambda_n^\rho)),\\ (2\lambda_n^{\rho})^
{\psi(2\lambda_n^{\rho})} & \leq &
\big(2\lambda_n^{\frac{\rho}{1-\chi(c^{-j})}}\big)^
{\psi\big(2\lambda_n^{\frac{\rho}{1-\chi(c^{-j})}}\big)}
\end{eqnarray*}
for some constant $M\ge 1$ also independent of $y$ and $j$. Eventually,
since $\chi(r) \ra 0$ when $r\ra 0$, one has $\frac{1}{1-\chi(c^{-j})}
\leq 1+ 2 \chi(c^{-j})$ for $j$ large enough. As a consequence, for
the same constant $M$ one can write
\[
M^{-1}(2\lambda_n^\rho)^{\alpha+2\alpha\chi(2\lambda_n^{\rho})+
\psi(2\lambda_n^{\rho})} \leq \mu(B(x_n, \lambda_n^\rho)) .\]

The upper bound of (\ref{wizzz}) is treated with the same arguments,
and one obtains $\displaystyle \mu(B(x_n, \lambda_n^\rho)) \leq M
(2\lambda_n^\rho)^{\alpha-\alpha \chi(2\lambda_n^{\rho})-
\psi(2\lambda_n^{\rho})}$. Hence $\mathcal{P}^\rho_{M}(B(x_n,
\lambda_n^\rho))$ holds.

\smallskip

To prove that $\mathcal{D}^m_{M}(B(y, r)$ holds for some $M>0$
independent of $y$ and $r\in (0, c^{-j(y)-1})$ it is enough to write
that $ B(y, r)\subset \bigcup_{\|{\bf k} - {\bf
k}_{j,y}^c\|_{\infty} \leq 1} I_{j,{\bf k}}^c$, where $j$ is the
largest integer such that $r\le c^{-j}$, and then to use (\ref{P1}).
\eproof

If $y$, $j$ and $(x_n,\lambda_n)$ satisfy (\ref{P0r}),
then they also satisfy~(\ref{wizr}). This ensures that
the Cantor set we are going to build is included in
$\smr$.

\bproof[of Theorem \ref{theor_lower} in the case $\rho<1$] Here again,
 the case $\delta=1$ is obvious and left to the reader. Since
 $\mathcal{D}= (1,\infty)$, we deal with the sets $\widehat
 S_\mu(\rho, (\delta)_{n\ge 1},\alpha,\varepsilon^\rho_M)$, which are
 equal to the sets $S_\mu(\rho, \delta,\alpha,\varepsilon^\rho_M)$.
\smallskip

Let $\delta>1$. As in the proof of Theorem \ref{theor_lower}, we
 construct a generalized Cantor set $K_\delta$ in $\smr$ and a measure
 $m_{\rho,\delta}$ on $K_\delta$.

\medskip

{\bf - First step:} The first generation in the construction of
$K_\delta$ is as~follows:

Let $L_0=[0,1]^d$. Using assumption {\bf (2($\rho$))}, there exist a
subset $E^{L_0}$ of $E^{L_0}_{n_{L_0}}$ of $m$-measure larger than
$\Vert m\Vert/4$ and an integer $n'_{L_0}\ge n_{L_0}$ such that
$\forall y\in E^{L_0}$, $\forall j\ge n'_{L_0}$, (\ref{P2}) and
(\ref{P1}) hold. There is a subset $\widetilde E^{L_0} $ of $ E^{L_0}$
of $m$-measure greater than $\Vert m\Vert/8$ such that for every $y
\in \widetilde E^{L_0}$, (\ref{P0r}) holds.

Once again we are going to apply Theorem \ref{besic} to $A=\widetilde
E^{L_0}$ and to families $\mathcal{B}_1(j)$ of balls built as follows.
Let $y \in \widetilde E^{L_0}$. We define
\begin{equation}
\label{njyr}
n_{j,y,\rho}= \inf \left\{n:\! c^{-n(1-\chi(c^{-n}))} \leq
c^{-\frac{j+5}{\rho}}\mbox{ and }\! (\ref{P0r}) \mbox{ holds 
with } j_i(y)=n \right\}.
\end{equation}
Then for every $j\geq n'_{L_0}$, let us introduce the family
\[{\mathcal{B}}_1(j)=\left\{B(y, 3c^{-\rho n_{j,y,\rho}}):\
y\in \widetilde E^{L_0}\right\}.
\]

For every $j\geq n'_{L_0}$, the family ${\mathcal{B}}_1(j)$ fulfills
conditions of Theorem \ref{besic}. 

Hence, $\forall j\geq n'_{L_0}$, $Q(d)$ families of disjoint balls
${\mathcal{B}}^1_1(j),...,{\mathcal{B}}_1^{Q(d)}(j)$ can be extracted
from ${\mathcal{B}}_1(j)$.  The same procedure as in Theorem
\ref{theor_lower} allows us to extract from these new families a
finite family of disjoint balls $\widetilde{G}_1(j) =
\{B_1,B_2,\ldots, B_N\}$ such that
\begin{equation}
\label{minormass1r}
m\Big( \bigcup _{B_k \in \widetilde{G}_1(j)} B_k \Big) \geq
\frac{\Vert m\Vert}{ 16 \,Q(d)}.
\end{equation}

Remember that with each $B_k$ can be associated a point $y_k\in
\widetilde E^{L_0} $ so that $B_k = B(y_k,3c^{-{\rho
n_{j,y_k,\rho}}})$.  Let us fix one of the balls $B_k =
B(y_k,3c^{-\rho n_{j,y_k,\rho}})$. By construction, one can find
$[c^{n_{j,y_k,\rho}(d(1-\rho)-\chi(c^{-n_{j,y_k,\rho}}))}]$ points
$x_n$ in the ball $B(y_k,c^{-\rho n_{j,y_k,\rho}})$ such that
(\ref{P0r}) holds. We denote $\mathcal{S}(B_k)$ the set of these
points $x_n$. The corresponding balls $B(x_n, \lambda_n)$ are pairwise
disjoint. By construction, for each of these points $x_n \in
\mathcal{S}(B_k)$, one has
\begin{equation}
\label{inegr}
B(y_k,(c^\rho-1)c^{-\rho n_{j,y_k,\rho}}) \subset B(x_n,
\lambda_n^\rho) \subset B\big(y_k, c^{-\rho n_{j,y_k,\rho}
(1-\chi(c^{-n_{j,y_k,\rho}}))}\big).
\end{equation}

Therefore each point $x_n\in \mathcal{S}(B_k)$ such that (\ref{P0r})
holds verifies the conditions of Lemma \ref{btolr}. Thus ${\mathcal
P}^\rho_{M}(B(x_n, \lambda_n^\rho))$ and $\mathcal{D}^m_M(B_k)$ hold
for some constant $M$ independent of the scale and of $x$. This
constant $M$ is the one chosen to define $\smr$.

Let us now consider $I_n^{(\delta)} = B(x_n, \lambda_n^\delta)$. Let
$J_{n,k}$ be the closure of one of the $c$-adic box of maximal
diameter included in $I_n^{(\delta)}$. Since $|B_k|= 6c^{-\rho
n_{j,y_k,\rho}}$, one has $ |B_k|\leq C|J_{n,k}|^{{\rho}/{\delta}}
$ for some constant $C$ depending only on $\delta$.

We write $\underline{B_k} = J_{n,k}$.  Conversely, if a closed
 $c$-adic box $J $ can be written $\underline{B} $ for some larger
 ball $B$, one writes $B=\overline J$. Pay attention to the fact that
 a number equal to $\#\mathcal{S}(B_k)\ge [
 c^{n_{j,y_k,\rho} (d(1-\rho)-\chi(c^{-n_{j,y_k,\rho}}))}]$ of $c$-adic
 boxes $J_{n,k}$ can be written as $\underline{B_k}$ for the same ball
 $B_k$.  For every $c$-adic box $J$ such that there exists $k$ with
 $B_k =\overline J$, one ensured by construction
\begin{equation}
\label{toto1r}
 |\overline J| \leq C|J|^{{\rho}/{\delta}}
\end{equation}
 for some constant $C$ depending on $\delta$. Moreover, the $c$-adic
 box $J$ is included in a contracted ball
 $I_n^{(\delta)}=B(x_n,\lambda^\delta_n)$ such that ${\mathcal
 P}^\rho_{M}(B(x_n,\lambda^\rho_n))$ holds.

Since $|B_k | = 6 c^{-\rho n_{j,y_k,\rho}} $, there is $C>0$
independent of $k$ and~$\rho$ such~that
\begin{equation}
\label{majs}
\#\mathcal{S}(B_k)\ge [ c^{n_{j,y_k,\rho}
 (d(1-\rho)-\chi(c^{-n_{j,y_k,\rho}}))} ]\geq {C}^{-1}
 |B_k|^{-\frac{d(1-\rho)}{\rho}} |B_k|^{\chi(|B_k|)}.
\end{equation}
We eventually define
\begin{equation}
\label{defg1jr}
G_1(j) = \{J_{n,k}: \overline{J_{n,k}} \in \widetilde{G}_1(j)\}.
\end{equation}
We notice that $I_1$ and $I_2$ belong to $G_1(j)$ and
$\overline{I_1}\neq \overline{I_2}$ then the distance between $I_1$
and $I_2$ is by construction at least $\max_{i\in
\{1,2\}} {\overline{I_i}}/{3}$. 

\medskip

On the algebra generated by the elements of $G_1(j)$, a probability
measure $m_{\delta,\rho}$ is defined by
\[ m_{\rho,\delta} (I)=\frac{\frac{m(\overline I)}{\#\, {\mathcal
S}(\overline I)}}{\sum_{B_k\in \widetilde G_1(j)}m( B_k)}.
\]
Since ${\mathcal D}_M^m(\overline I)$ holds for the measure $m$,
by (\ref{toto1r}) and (\ref{monotone}), we have
\[
 m(\overline I) \le M |\overline I|^{\beta-\varphi (|\overline I|)}\le
C|I|^{{\rho\beta}/{\delta}}|\overline I|^{-\varphi(|\overline I|)}
\le C|I|^{{\rho\beta}/{\delta}}|I|^{-\varphi(| I|)}.
\] Then, one also has by (\ref{majs}) and (\ref{inegr})
\begin{eqnarray*}
\label{maj2r}
 ({\#\, {\mathcal S}(\overline I)})^{-1}& \le & C|\overline
I|^{\frac{d(1-\rho)}{\rho}}|\overline I|^{-\chi(|\overline I|)} \le C
| I|^{\frac{\rho}{\delta}\frac{d(1-\rho)}{\rho}}| I|^{-\chi(| I|)}
\leq C | I|^{\frac{d(1-\rho)}{\delta}}| I|^{-\chi(| I|)}.
\end{eqnarray*}
Moreover, by
(\ref{minormass1r}) and the definition of $G_1(j)$ (\ref{defg1j}), one
gets
\[
{\sum_{B_k\in \widetilde G_1(j)} m(B_k)} \ge \frac{\Vert m\Vert}{16\,
Q(d)}.
\]
Thus, $\displaystyle \forall \ I\in G_1(j), \,\,m_{\rho,\delta} (I)\le
{16\, Q(d)C{\Vert m\Vert}^{-1}| I|^{-\varphi(|I|)}| I|^{-\chi(| I|)}}
|I|^{\frac{d(1-\rho)+\rho \beta}{\delta}}$.  By our assumption {\bf
(1)}, we can fix $j_1$ large enough so that
\[
\forall \ I\in G_1(j_1), \ {16 \,Q(d)C}{\Vert m\Vert}^{-1} \le |
I|^{-\varphi(|I|)} .
\] We choose the $c$-adic elements of the first generation of the
construction of $K_\delta$ as being those of $G_1:=G_1(j_1)$. By
construction
\begin{equation}
\label{premmajr} \forall \ I\in G_1, \, \, m_{\rho,\delta} (I)\le
|I|^{\frac{d(1-\rho)+\rho\beta}{\delta}-2\varphi(|I|)- \chi(| I|)},
\end{equation}
and for every $x\in \bigcup_{I\in G_1}I$, there exists an integer $n$
so that $ \lambda_n \leq c^{-{5}/{\rho}} $, \\$\|x_n-x\|_\infty\leq
\lambda_n^\delta$, and ${\mathcal P}^\rho_{M}(B(x_n, \lambda_n^\rho))$
holds. Moreover, $\max_{I\in G_1}|I|\le 2c^{-{5\delta}/{\rho}}$.

\medskip

{\bf  -  Second step:} 
The second generation is built as in the case $\rho=1$, by focusing on
one $c$-adic box $L$ of the first generation. We give the essential clues
to obtain this second generation.

Using assumption {\bf (2($\rho$))}, there exist a subset $E^{L}$ of
$E^{L}_{n_{L}}$ of $m^L$-measure larger than $\Vert
m^L\Vert/4$ and an integer $n'_{L}\ge n_{L}$ such that for all $y\in
E^{L}$, for every $j\ge n'_{L}+\log_c\big (|L|^{-1}\big )$,
(\ref{ref4}) holds. Then, there exists a subset $\widetilde E^{L} $ of
$ E^{L}$ of $m^L$-measure greater than $\Vert m^L\Vert/8$
such that for every $y \in
\widetilde E^{L}$, (\ref{P0r}) holds.

One more time we apply Theorem \ref{besic} to $A=\widetilde E^{L}$ and
to families of balls $\mathcal{B}_2(j)$.  Let $y \in \widetilde
E^{L}$. For every $j\geq n'_{L}+\log_c\big (|L|^{-1}\big )$, we define
the family
\[{\mathcal{B}}_2(j)=\Big\{B(y, 3c^{-\rho n_{j,y,\rho}}):\ y\in
\widetilde E^{L}\Big\}.
\]
The family $\widetilde{\mathcal{B}}_2(j)$ fulfills conditions of
Theorem \ref{besic}. Hence, $Q(d)$ families of disjoint balls
${\mathcal{B}}^1_2(j),...,{\mathcal{B}}_2^{Q(d)}(j)$ can be extracted
from ${\mathcal{B}}_2(j)$. Moreover, one can also extract from these
families one finite family of disjoint balls $\widetilde{G}^L_2(j) =
\{B_1,B_2,\ldots, B_N\}$ such that
\begin{equation}
\label{minormass2r}
m^L\Big( \bigcup _{B_k \in \widetilde{G}_2(j)} B_k \Big) \geq
\frac{\Vert m^L\Vert}{ 16 \, Q(d)}.
\end{equation}

Each of these balls $B_k$ can be written $B(y_k,3c^{-\rho
n_{j,y_k,\rho}})$ for some point $y_k\in\widetilde E^L$ and some
integer $n_{j,y_k,\rho}$. Moreover, by (\ref{P0r}), with each $B_k$
can be associated
$[c^{n_{j,y_k,\rho}(d(1-\rho)-\chi(c^{-n_{j,y_k,\rho}}))}]$ points
$x_n$ in $B(y_k,c^{-\rho n_{j,y_k,\rho}})$ such that (\ref{P0r})
holds. As above, $\mathcal{S}(B_k)$ denotes the set of these points
$x_n$. The corresponding balls $B(x_n, \lambda_n)$ are pairwise
disjoint.

By construction, (\ref{inegr}) holds for each of these points $x_n \in
 \mathcal{S}(B_k)$. Moreover, Lemma \ref{btolr} holds with the measure
 $m^L\circ f_L^{-1}$ instead of $m$ and with the same constant
 $M$. Consequently, each point $x_n\in \mathcal{S}(B_k)$ such that
 (\ref{P0r}) holds is such that ${\mathcal P}^\rho_{M}(B(x_n,
 \lambda_n^\rho))$ and $\mathcal{D}^{m^L\circ f_L^{-1}}_M\big
 (f_L(B_k)\big )$ hold.

We then consider $I_n^{(\delta)} = B(x_n, \lambda_n^\delta)$, and we
denote by $J_{n,k}$ the closure of one $c$-adic box of maximal
diameter included in $I_n^{(\delta)}$. Again one has (\ref{toto1r}).

We write $\underline{B_k} = J_{n,k}$.  Conversely, if a closed
 $c$-adic box $J $ can be written $\underline{B} $ for some larger
 ball $B$, one writes $B=\overline J$.  We eventually set
\begin{equation}
\label{defg2jr}
G^L_2(j) = \{J_{n,k}: \overline{J_{n,k}} \in \widetilde{G}^L_2(j)\}.
\end{equation}
 
On the algebra generated by the elements of $G_2^L(j)$, an extension
of the probability measure $m_{\rho,\delta}$ is defined by
\[ m_{\rho,\delta} (I)=m_{\rho,\delta} (L)\frac{ \frac{m^L(\overline
I)}{\# {\mathcal S}(\overline I)}}{\sum_{B_k \in \widetilde
G^L_2(j)}m^L(B_k)}.
\]
Since $\mathcal{D}^{m^L\circ f_L^{-1}}_M\big (f_L(B_k)\big )$ and
(\ref{toto1r}) hold, one gets
\begin{eqnarray*}
m^L(\overline I) \leq \left (\frac{|\overline I|}{|L|}\right
)^{\beta-\varphi \big(\frac{|\overline I|}{|L|}\big)} \!\!  \!\! \leq
C|I|^{\frac{\rho \beta}{\delta}}|L|^{-\beta}\left (\frac{|\overline
I|}{|L|}\right )^{-\varphi \big(\frac{|\overline I|}{|L|}\big)} \!\!
\!\! \leq C|I|^{\frac{\rho \beta}{\delta}}|L|^{-\beta} |I|^{-\varphi
(|I|)},
\end{eqnarray*}
where the monotonicity of $x\mapsto x^{-\varphi(x)}$ of assumption
{\bf (1)} is used. Then (\ref{majs}) applied to $\overline{I}$ and
(\ref{minormass2r}) yield
\begin{eqnarray*}
 m_{\rho,\delta} (I)&\le & m_{\rho,\delta}(L) \frac{16\,Q(d) C}{\Vert
m^L\Vert} |I|^{\frac{\rho\beta}{\delta}}|L|^{-\beta}
|I|^{-\varphi (|I|)}
|I|^{\frac{d(1-\rho)}{\delta}}| I|^{-\chi(| I|)},
\end{eqnarray*}
and using (\ref{premmajr}) finally gives
\begin{eqnarray*}
 m_{\rho,\delta} (I)&\le& \frac{16\,Q(d) C |L|^
{\frac{d(1-\rho)+\rho\beta}{\delta}-\beta-2\varphi(|L|)- \chi(|
L|)}}{\Vert m^L\Vert}
|I|^{\frac{d(1-\rho)+\rho\beta}{\delta}-\varphi (|I|)-\chi(| I|)}
\end{eqnarray*}
By assumption {\bf (1)} one can choose $j_2(L)$ large enough so that
for every integer $j\ge j_2(L)$, for every $I\in G_2^L(j)$,
\[
{16\,Q(d) C{\Vert
m^L\Vert}^{-1}|L|^{\frac{d(1-\rho)+\rho\beta}{\delta}-\beta-2\varphi
(|L|)- \chi(|L|)}} \le |I|^{-\varphi (|I|)}.
\]
Then, taking $j_2=\max\big \{j_2(L):\ L\in G_1\big\}$ and defining $
G_2=\bigcup_{L\in G_1}G^L_2(j_2), $ this yields an extension of
$m_{\rho,\delta}$ to the algebra generated by the elements of
$G_1\bigcup G_2$. One has for every $I\in G_1\bigcup G_2$, $
m_{\rho,\delta}(I)\le |I|^{\frac{d(1-\rho)+\rho\beta}{\delta}-2\varphi
(|I|)-\chi(| I|)}.  $

Remark that by construction if $J\in G_1$ and $I\in G_{2}$ verify
$I\subset J$ one has
\[ \sum_{I'\in G_{2},\ \overline I'=\overline {I}} m_{\rho,\delta}
(I')\le 16\,Q(d) m_{\rho,\delta}(J) \frac{ m^{J}(\overline I)}{\Vert
m^{J}\Vert}.
\]

Also notice that by construction, $|I|\le \max_{J\in G_1}
2(c^{-5}|J|)^{{\delta}/{\rho}}\le (2c^{-5{\delta}/{\rho}})^2$
for every $I\in G_2$. Moreover, $I$ is contained in some
$I_n^{(\delta)}$ such that $|I_n^{(\delta)}|\le C |I|$, where $C$ is a
constant which depends only on $c$.

\medskip

{\bf - Third step:} Assume that $n$ generations of closed $c$-adic
boxes $G_1,\dots, G_n$ have already been found for some integer $n\ge
2$. Assume also that a probability measure $m_{\rho,\delta}$ on the
algebra generated by $\bigcup_{1\le p\le n}G_p$ is defined and that:

\smallskip

{\bf (i)} The elements of $G_p$ are pairwise disjoint closed $c$-adic
boxes, and for $1\le p\le n$, $\max_{I\in G_p}|I|\le \big
(2c^{-5{\delta}/{\rho}}\big )^p$.

For every $1\le p\le n$, with each $I\in G_p$ is associated a ball
$\overline{I}$ such that $I\subset \overline{I}$. There exists a
constant $C>0$ which depends only on $\delta$ such that (\ref{toto1r})
holds. Moreover, if $I_1$ and $I_2$ belong to $G_p$ and
$\overline{I_1}\neq\overline{I_2}$, their distance is at least
$\max_{i\in\{1,2\}}{\overline I_i}/{3}$. Moreover, the
$\overline{I}$'s ($I\in G_p$) are pairwise disjoint.

\smallskip

{\bf (ii)} For every $2\le p\le n$, each element $I$ of $G_p$ is a
subset of an element $L$ of $G_{p-1}$. Moreover, $\overline I\subset
L$, $\log_c \big (|\overline I|^{-1}\big )\ge n_{L}+\log_c \big
(|L|^{-1}\big )$ and $\overline I\cap E^{L}_{n_{L }}\neq\emptyset$.

\smallskip

{\bf (iii)} For every $1\le p\le n$ and $I\in G_p$, there exists an
integer $q$ such that $I\subset B(x_q,\lambda_q^\delta)=
I_q^{(\delta)}\subset \overline{I}$ and ${\mathcal
P}^\rho_{M}(B(x_q,\lambda_q^\rho))$ holds, and $|I_q^{(\delta)}|\le
C|I|$ for some constant $C$ which depends only on $c$.

{\bf (iv)} For every $I\in \bigcup_{1\le p \le n}G_p$, 
$
m_{\rho,\delta}(I)\le
|I|^{\frac{d(1-\rho)+\rho\beta}{\delta}-2\varphi(|I|)- \chi(|I|)}.$

\smallskip

{\bf (v)} For every $1\le p\le n-1$, $L\in G_p$, and $I\in G_{p+1}$ such
that $I\subset L$,
\[ \sum_{I'\in G_{p+1},\ \overline{I'}=\overline {I}} m_{\rho,\delta}
(I')\le 16 \, Q(d) m_{\rho,\delta}(L) \frac{ m^{L}(\overline I)}{\Vert
m^{L}\Vert}.
\]
The construction of a generation $G_{n+1}$ of $c$-adic boxes and an
extension of $m_{\rho,\delta}$ to the algebra generated by the
elements of $\bigcup_{1\le p\le n+1} G_p$ such that properties {\bf
(i)} to {\bf (v)} hold for $n+1$ are done as when $n=1$. 

Then, by induction, we get a sequence $(G_n)_{n\ge 1}$ and a
probability measure on $\sigma \big(I:\ I\in \bigcup_{n\ge 1}G_j\big
)$ such that properties {\bf (i)} to {\bf (v)} hold for every $n\ge
2$, and $\displaystyle K_{\rho,\delta}=\bigcap_{n\ge 1} \,
\bigcup_{I\in G_n}I$. By construction,
$m_{\rho,\delta}(K_{\rho,\delta})=1$ and because of {\bf (iii)}
$K_{\rho,\delta}\subset \smr$. Finally, the measure $m_{\rho,\delta}$
is extended to $\mathcal{B}([0,1]^d)$ in the usual way:
$m_{\rho,\delta}(B):=m_{\rho,\delta}(B\cap K_{\rho,\delta})$ for every
$B\in \mathcal{B}([0,1]^d)$.

\medskip

{\bf - Last step:} Proof of (\ref{controldemdelta}). If $I \in
G_n$, recall that we set $g(I) =n $.

Fix $B$ an open ball of $[0,1]$ of diameter less than the one of the
elements of $G_1$ such that $B\cap K_{\rho,\delta} \neq
\emptyset$. Let $L$ be the element of largest diameter in
$\bigcup_{n\ge 1}G_n$ such that $B$ intersects at least two balls
$\overline L_i$ such that $L_i$ belongs to $G_{g(L)+1}$ and $ L_i$ is
included in $L$ (hence $ m_{\rho,\delta} (B)\le m_{\rho,\delta} (L)$).

\smallskip

\noindent
$\bullet$ If $|B|\ge |L|$,
\begin{eqnarray*}
 m_{\rho,\delta }(B) \!  \le \!  m_{\rho,\delta} (L) \!  \le \! 
|L|^{\frac{d(1-\rho)+\rho\beta}{\delta}-2\varphi (|L|)- \chi(|L|)} \!  \leq \!  C |B|^{\frac{d(1-\rho)+\rho\beta}{\delta}-2\varphi (|B|)-
   \chi(|B|)}.
\end{eqnarray*}

\smallskip
\noindent
$\bullet$ If $|B|<c^{-n_{L}-3}|L|$, let $L_1,\dots,L_p$ be the
$c$-adic boxes in $G_{g(L)+1}$ such that $\forall i$ $\overline L_i$
intersects $B$. Property {\bf (v)} yields
\begin{eqnarray*}
m_{\rho,\delta}(B) = \sum_{i=1}^p \sum_{L\in G_{g(L)+1},\ \overline
L=\overline L_i}m_{\rho,\delta} (B\cap L) \leq \sum_{i=1}^p m_{\rho,\delta}(L)
\frac{16\, Q(d)}{\Vert m^L\Vert} m^{L}(\overline L_i).
\end{eqnarray*}
Let $j_0$ be the unique integer so that $c^{-j_0}\le |B| <
c^{-j_0+1}$. Because of {\bf (i)}, one has $|B|\ge {\max_i
|\overline{L}_i|}/{3}$.  As a consequence $-\log_c |\overline{L}_i|
\geq j_0 -[\log_c 3] \geq j_0 -2$.

\smallskip

The same arguments as in the proof of Theorem \ref{theor_lower} (Case
$\rho =1$) yield that there exists an index $i_0$ and a point $y\in
E^L_{n_L}\cap \overline L_{i_0}$ such that one has 
$\bigcup_{i=1}^p \overline L_{i} \subset \bigcup_{{\bf
k}:\ \|{\bf k}-{\bf k}_{j_0-3,y}\|_{\infty} \leq 1} I_{j_0-3,{\bf
k}}^c$. Hence
\begin{equation}
\label{eq2}
\sum_{i=1}^p m^{L}(\overline L_i) \le \sum_{{\bf k}:\ \|{\bf k}-{\bf
k}_{j_0-3,y}\|_{\infty} \leq 1} m^{L}(I^c_{j_0-3,{\bf k}}),
\end{equation}
and by definition of $E^L_{n_L}$, one can bound $m^{L}(I^c_{j_0-3,{\bf
k}})$ by
\[ m^{L}(I^c_{j_0-3,{\bf k}})\le \left (\frac{|I^c_{j_0-3,{\bf
k}}|}{|L|}\right )^{\beta-\varphi\big(\frac{|I^c_{j_0-3,{\bf
k}}|}{|L|}\big)}\le C\left (\frac{|B|}{|L|}\right )^{\beta}\left
(\frac{|B|}{|L|}\right )^{-\varphi\big(\frac{|B|}{|L|}\big)}.
\] 
There are $3^d$ such pairwise disjoint boxes in the sum (\ref{eq2}),
hence
\begin{eqnarray*}
m_{\rho,\delta} (B)&\le& \frac{16\, Q(d) }{\Vert
m^{L}\Vert}m_{\rho,\delta} (L) 3^d C\left (\frac{|B|}{|L|}\right
)^{\beta}\left (\frac{|B|}{|L|}\right
)^{-\varphi\big(\frac{|B|}{|L|}\big)}\\ &\le& \frac{16\, Q(d) 3^d
C}{\Vert m^{L}\Vert}m_{\rho,\delta} (L) \left (\frac{|B|}{|L|}\right
)^{\beta}|B|^{-\varphi (|B|)}.
\end{eqnarray*}
By {\bf (iv)}, one obtains
\begin{eqnarray*}
 m_{\rho,\delta} (L)& \le & |L|^{\frac{d(1-\rho)+\rho\beta}{\delta}}
 |L|^{-2\varphi(|L|)-\chi(|L|)}\leq 
 |L|^{\frac{d(1-\rho)+\rho\beta}{\delta}}
 |B|^{-2\varphi(|B|)-\chi(|B|)},
\end{eqnarray*}
which yields
\[ m_{\rho,\delta} (B)\le \frac{16\, Q(d) 3^d C}{\Vert
m^L\Vert}|L|^{\frac{d(1-\rho)+\rho\beta}{\delta}} \left
(\frac{|B|}{|L|}\right )^{\beta}
|B|^{-3\varphi(|B|)- \chi(|B|)}.
\]
Then, the second property of (\ref{P3}) in assumption {\bf (4)} allows
to upper bound ${\Vert m^L\Vert}^{-1}$ by $|L|^{-\varphi(|L|)}$, which
is lower than $|B|^{-\varphi(|B|)}$, and thus
\begin{equation}\label{cas1rho}
m_{\rho,\delta} (B)\le C|L|^{\frac{d(1-\rho)+\rho\beta}{\delta}} \left
(\frac{|B|}{|L|}\right )^{\beta}|B|^{-4\varphi(|B|)
- \chi(|B|)}.
\end{equation}
Finally, if $\beta>\frac{d(1-\rho)+\rho\beta}{\delta}$, (\ref{cas1rho})
yields
\begin{eqnarray*}
m_{\rho,\delta} (B) & \le & C|B|^{\frac{d(1-\rho)+\rho\beta}{\delta}}
\left (\frac{|B|}{|L|}\right )^{\beta-\frac{d(1-\rho)+\rho\beta}{\delta}}
|B|^{-4\varphi(|B|)- \chi(|B|)}\\
&\le & C|B|^{\frac{d(1-\rho)+\rho\beta}{\delta}} |B|^{-4\varphi(|B|)-
  \chi(|B|)};
\end{eqnarray*}
If $\beta\le \frac{d(1-\rho)+\rho\beta}{\delta}$, (\ref{cas1rho}) yields
\begin{eqnarray*}
 m_{\rho,\delta} (B)& \le & C |B|^\beta
 |L|^{\frac{d(1-\rho)+\rho\beta}{\delta}-\beta}
 |B|^{-4\varphi(|B|)- \chi(|B|)} \leq C |B|^\beta |B|^{-4\varphi(|B|)- \chi(|B|)}.
\end{eqnarray*}
In both cases, if $D(\beta,\rho,\delta)=\min(\beta,
\frac{1-\rho+\rho\beta}{\delta})$,
\begin{equation}
\label{finalr}
m_{\rho,\delta} (B)\le C |B|^{D(\beta,\rho,\delta)}
|B|^{-4\varphi(|B|) -\chi(|B|)}.
\end{equation}

$\bullet$ $c^{-n_L-3}|L|\le |B|\le |L|$: one needs at most $c^{d(n_L+4)}$
contiguous $c$-adic boxes of length $c^{-n_L-3}|L|$ to cover $B$. For these
boxes, (\ref{finalr}) can be used to~get
\begin{eqnarray*}
m_{\rho,\delta} (B)&\le& C c^{d(n_L+4)} \big (c^{-n_L-3}|L|\big
)^{D(\beta,\rho,\delta)-4\varphi(c^{-n_L-3}|L|)-\chi(c^{-n_L-3}|L|)}\\
&\le &C c^{dn_L} |B|^{D(\beta,\rho,\delta)}
|B|^{-4\varphi(|B|)-\chi(|B|)} \\
&\le & C|L|^{-d\varphi (|L|)}|B|^{D(\beta,\rho,\delta)}
|B|^{-4\varphi(|B|)-\chi(|B|)}\\
&\le &C |B|^{D(\beta,\rho,\delta)} |B|^{-(4+d)\varphi(|B|)-\chi(|B|)}.
\end{eqnarray*}
This shows (\ref{controldemdelta}) and ends the proof of Theorem
\ref{theor_lower} when $\rho<1$.  \eproof

\section{Examples}
\label{sec_examples}

Section~\ref{sec_couples} exhibits several families
 $\{(x_n,\lambda_n)\}_{n}$ which satisfy (\ref{P0}) or (\ref{P0r}) for
 any measure $m$, and form weakly redundant systems. Then
 Section~\ref{sec_measures} provides examples of triplets $\big
 (\mu,\alpha,\tau_\mu^*(\alpha)\big )$ leading to $\rho$-heterogeneous
 ubiquitous systems. It also gives relevant interpretations to
 property~$\mathcal{P}^\rho_M$.

\subsection{Examples of families $\{(x_n,\lambda_n)\}_{n\in \mathbb{N}}$}
\label{sec_couples}

 Let us notice first that, to ensure (\ref{P0}), it suffices that
\begin{equation}
\label{recouvre1}
\bigcap_{N \geq 1} \bigcup_{n\geq N} B\big(x_n,{\lambda_n}/{2}\big)
=[0,1]^d.
\end{equation} 

 {$\bullet$ Family of the $b$-adic numbers. } 

\medskip Fix $b$ an integer $\ge 2$. Let us consider the sequence
$\{({\bf k}b^{-j},2b^{-j})\}$, for $j\in\mathbb{N}$ and ${\bf k}=(k_1,
k_2,\ldots, k_d)\in \{0,\ldots,b^{j}-1\}^d$. By construction, for
every $j\geq 2$, $\bigcup_{{\bf k}\in \{0,\ldots,b^j-1\}^d} B\big({\bf
k}b^{-j}, b^{-j}\big) = [0,1]^d$. Hence (\ref{recouvre1}) is
satisfied, (\ref{P0r}) holds for any measure $m$ and the family is
weakly redundant.

%

\medskip 

{$\bullet$ Family of the rational numbers.}

\medskip By Theorem 200 of \cite{HARDYWRIGHT}, any point
$x=(x_1,\ldots,x_d)\in [0,1]^d$ such that at least one of the $x_i$ is
an irrational number satisfies for infinitely many ${\bf p}=(p_1,
p_2,\ldots, p_d)$ and $q$ the inequality $\|x-{{\bf p}}/{q}\|_{\infty}
\leq {q^{-(1+{1}/{d})}}$. As a consequence, the sequence
$\big\{\big({{\bf p}}/{q}, {2}{q^{-(1+{1}/{d})}}\big)\big\}$ for $
q\in \mathbb{N^*}$ and ${\bf p}=(p_1, p_2,\ldots, p_d)\in
\{0,\ldots,q-1\}^d$ fulfills (\ref{recouvre1}). Here again,
(\ref{P0r}) holds for any measure $m$.

\smallskip

To ensure the weak redundancy, one must select only the rational
numbers $\big\{\big({{\bf p}}/{q}, {2}{q^{-(1+{1}/{d})}}\big)\big\}$
such that at least one fraction ${p_i}/{q}$ is irreducible. But
(\ref{recouvre1}) is no more satisfied. Indeed, the rational numbers
${{\bf p}}/{q}$ themselves do not belong to the corresponding
limsup-set (each rational number belongs only to a finite number of
balls $B\big ({{\bf p}}/{q},
{2}{q^{-(1+{1}/{d})}}\big)$. Nevertheless, as soon as the rational
points are not atoms of $m$ (for instance if $\underline{\dim}(m)>0$),
both (\ref{P0}) and (\ref{P0r}) hold. In this case, by Theorem 193 of
\cite{HARDYWRIGHT}, the same holds with $\big\{\big({{ p}}/{q},
{2}/{\sqrt{5}q^{2}}\big)\big\}$ when $d=1$. This family is used to
prove (\ref{JaBe}).

\medskip 

 {$\bullet$ Family of the $\big \{(\{n\alpha\}, 1/n)\big\}_{n\in
    \mathbb{N}}$.}

\medskip
 Let us focus on the case $d=1$ to introduce another
family. Let $\alpha$ be an irrational number. For every $n\in
\mathbb{N}$, we denote by $\{n \alpha\}$ the fractional part of
$n\alpha$. If $x \notin \mathbb{Z}+\alpha\mathbb{Z}$, one has
$|n\alpha-x| < {1}/{2n}$ for an infinite number of integers $n$
(see Theorem II.B in \cite{cassels} for instance). Hence
\[\mathbb{R} \bs\left(\mathbb{Z}+\alpha\mathbb{Z}\right) \subset
\bigcap_{N\geq 1} \bigcup_{n\geq N} B(\{n\alpha\},{1}/{2n}).\]
As soon as $m\left(\mathbb{Z}+\alpha\mathbb{Z}\right) =0$, (\ref{P0})
is satisfied for the family $\{(\{n\alpha\},1/n)\}_{n\ge 1}$.
We do not know the measures $m$ for which (\ref{P0r}) holds.  However the
following property concerning the redundancy holds:
\begin{proposition}
$\{(\{n\alpha\},1/n)\}_{n\ge 1}$ forms a weakly redundant system if
and only if $\ \displaystyle \inf\left\{\xi:\
\#\left\{(p,q)\in\mathbb{N}\times\mathbb{N}^*:\ \left
|\alpha-{p}/{q}\right |\le q^{-\xi}\right \}=\infty\right \}=2$.
\end{proposition}
One knows that every irrational number is approximated at rate
$\xi\geq 2$ by the rational numbers. But the system
$\{(\{n\alpha\},{1}/{n})\}_{n}$ is weakly redundant if and only if the
approximation rate by rational numbers of $\alpha$ is exactly
equals~2.

\begin{proof}
Notations of Definition~\ref{weaksystem} are used. 

Remark that $T_j$ (defined by (\ref{deftj})) contains exactly $2^j$
integers.

Suppose that the family is not weakly redundant. For every partition
of $T_j$ into $N_j$ subsets, one has $ \limsup_{j\ra +\infty }
j^{-1}{\log N_j} > 0$. Let us fix such a partition. There exists
$\ep>0$ such that for infinitely many integers $j$, one can find a
real number $x\in [0,1]$ such that more than $2^{ \ep j}$ among the
$\{B(x_n,\lambda_n)\}_{n\in T_j}$ contain $x$. Since these integers
$n$ belong to $T_j$, the corresponding $\lambda_n$ belong to
$(2^{-(j+1)},2^{-j}]$. Consequently, these $2^{ \ep j}$ integers $n$
all verify $|\{n\alpha\} -x |\leq 2^{-j}$.

By a classical argument, there are two integers $n$ and $n'$ of $T_j$
such that
\begin{equation}
\label{key}
n\neq n', \ |n-n'|\le 2^{j} \mbox{ and } |\{n\alpha\}-\{n'\alpha\}|\le
2\cdot2^{-j(1+\ep)}.
\end{equation}
We deduce from (\ref{key}) that there exists $p\in\mathbb{N}$ such
that $ \big ||n-n'|\,\alpha-p\big |\le 2\cdot 2^{-j(1+\ep)} \le
2|n-n'|^{-(1+\ep)}$. Hence $\big |\alpha-{p}/{|n-n'|}\big |\le
2|n-n'|^{-(2+\ep)}$. Since (\ref{key}) holds for infinitely many $j$,
$|n-n'|$ cannot be bounded as $j$ goes to $\infty$. This yields
$\xi_\alpha:= \inf\big\{\xi:\
\#\big\{(p,q)\in\mathbb{N}\times\mathbb{N}^*:\ \big
|\alpha-{p}/{q}\big |\le q^{-\xi}\big \}=\infty\big \}>2$.

\smallskip

Conversely, if $\xi_\alpha>2$, fix $\varepsilon\in
(0,\xi_\alpha-2)$. For infinitely many
$(p,q)\in\mathbb{N}\times\mathbb{N}^*$, one has $\left
|\alpha-{p}/{q}\right |\le {q^{-(2+\varepsilon)}}$. For such an
integer $q$, one has $\{nq\alpha\}\le {1}/{qn}$ for every $n\in
\left [1, q^{{\varepsilon}/{2}}\right ]$. For $q$ large enough,
let $j_q$ be the largest integer $j$ so that $[j,j+1]\subset \left
[\log_2(q),(1+{\varepsilon}/{2}) \log_2(q)\right ]$. Consider then
$T_{j_q}$. By construction, the point $0$ belongs to at least
$2^{\frac{\ep}{4} {j_q}}$ balls $B(x_n,\lambda_n)$ such that ${n\in
T_{j_q}}$. Hence $N_{j_q} \geq 2^{{j_q}{\ep}/{4} }$. Since this
holds for infinitely many $j$'s, the conclusion follows.
\end{proof}

\medskip 
{$\bullet$ Poisson point processes.} 

\medskip
Let $S$ be a Poisson point process with intensity $\lambda\otimes \nu$
in the square $[0,1]\times (0,1]$, where $\lambda$ denotes the
Lebesgue measure on $[0,1]$ and $\nu$ is a positive locally finite
Borel measure on $(0,1]$ (see \cite{Kingman} for the construction of a
Poisson process). Let us take the family $\{(x_n,\lambda_n)\}_n$ equal
to the set $S$.  Let $c$ be an integer $\ge 2$. Then for $j\ge 1$, let
us introduce the quantities $T_j^c=\{n:\ c^{-(j+1)}<\lambda_n\le
c^{-j}\}$, as well as
\[
\beta_j=j^{-1}{\log_c\nu((c^{-(j-1)},c^{-(j-2)}])} \ \mbox{ and } \
 \beta=\limsup_{j\to\infty}\beta_j.
\]
One has $ \beta= \limsup_{j\to\infty}j^{-1}{\log_b\mathbb{E}(\#\,
T_{j-2})}$ for $b\in\{2,c\}$, but we use a basis $c$ rather than 2 in
order to discuss property (\ref{P0r}).  In fact, it is a general
property that the number $\limsup_{j\to\infty}j^{-1}{\log_c\#\,
T^c_j}$ itself does not depend on $c$.
We group the information concerning (\ref{P0}), (\ref{P0r}) and weak
redundancy:
\begin{proposition}\label{bbb}
\begin{enumerate}
\item 
Suppose $\int_{[0,1]}\exp\Big (2\int_{[t,1]}\nu ((2y,1))\, dy\Big ) dt
=+\infty$. This implies in particular $\beta\ge 1$. With probability 1,
(\ref{recouvre1}) holds.

\smallskip

\item
Fix $\rho\in (0,1)$. Let $\chi$ be a function defined as in
Definition~\ref{deffubir}. If there exists an increasing sequence
$(j_n)_{n\ge 1}$ such that $\beta_{j_n}\ge
1-\chi(c^{-j_n})+{4}/{j_n}$, then with probability 1, (\ref{P0r})
holds for any measure $m$.

\smallskip

\item $\{(x_n,\lambda_n)\}_n$ is weakly redundant almost surely if and
only if $\beta \le 1$.
\end{enumerate}
\end{proposition}
As a consequence, if $\nu (d\lambda)=\gamma{d\lambda}/{\lambda^2}$
with $\gamma>{1}/{2}$, with probability 1, the system $S$ is weakly
redundant and (\ref{recouvre1}) holds. In addition, if $\gamma$ is
large enough, with probability 1, (\ref{P0r}) holds for any measure
$m$.

\begin{proof} (i) It is a consequence of Shepp's theorem (see
\cite{Shepp} and \cite{bertoin}).

\smallskip

\noindent
(ii) We shall need the following lemma.
\begin{lemma}\label{oui!}
Let $\gamma\in \big ({1},{2},1\big )$. Let $N$ be a Poisson random
variable with parameter $M$. For all $p\ge 1$, one has $
\mathbb{P}(N\le M-M^\gamma)=O(M^{-p})\ \ (M\to\infty)$. 
\end{lemma}
The proof of Lemma \ref{oui!} uses the identity $\sum_{k=0}^n \exp
(-M)\frac{M^k}{k!}=\int_{M}^\infty \frac{u^n}{n!}e^{-u}\, du$ ($M>0$,
$n\in \mathbb{N}$) as well as Laplace's method for equivalents of
integrals.

For $j\ge 1$ and $0\le k\le c^{[j\rho]}-1$, let $\widehat
I^c_{[j\rho],k}$ be the subset of $I^c_{[j\rho],k}$ obtained by
keeping one over $c$ of the consecutive $c$-adic subintervals of
$I_{[j\rho],k}$ of generation $j-2$, that is $\displaystyle \widehat
I^c_{[j\rho],k}=\bigcup_{k'=0,\ldots,c^{j-[j\rho]-3}-1} I^c_{j-2,
c^{j-2-[j\rho]}k+ ck'}$. Let us also define the random sets $
\displaystyle S_{j,k}=\big \{n:\ \lambda_n\in
(c^{-(j-1)},c^{-(j-2)}],\ x_n\in \widehat I^c_{[j\rho],k}\big \}$, and
the random variables $N_{j,k}=\#\, S_{j,k}$. The $N_{j,k}$'s are
mutually independent Poisson random variables with parameter $M_j$
equal to the product of $\nu \big ((c^{-(j-1)},c^{-(j-2)}]\big )$ with
$ \big |\widehat I^c_{[j\rho],k}\big |$, that is
$M_j=c^{j\beta_j}\cdot c^{-[j\rho]-1}$.

 Fix $\gamma\in ({1}/{2},1 )$ and let $\displaystyle E_j=\big
\{\forall\ 0\le k\le c^{[j\rho]}-1,\ N_{j,k}\ge M_j-M_j^\gamma\big \}$
for $j\ge 1$.  One has $\mathbb{P} (E_j)=\big (\mathbb{P}(N_{j,0}\ge
M_j-M_j^\gamma)\big ) ^{c^{[j\rho]}}$. Moreover, by definition of
$j_n$, one has $\lim_{n\to\infty} M_{j_n}=\infty$. Consequently, using
the form of $M_j$ and Lemma~\ref{oui!}, one has
$\lim_{n\to\infty}\mathbb{P} (E_{j_n})=1$. Since the events $E_{j_n}$
are independent, by the Borel-Cantelli lemma one has
$\mathbb{P}(\limsup_{n\to\infty}E_{j_n})=1$.

A computation shows that $M_{j_n}-M_{j_n}^\gamma\ge
c^{(\beta_{j_n}-\rho)j_n-4}$ for $n$ large enough. It follows that with
probability 1, there exist infinitely many $j_n$ such that for all
$0\le k\le c^{[j_n\rho]}-1$, $N_{j_n,k}\ge c^{j_n(1-\rho-\chi
(c^{-j_n}))}$. Moreover, by construction, the balls $B(x_n,\lambda_n)$
for $n\in S_{j,k}$ are pairwise disjoint, and if $y\in [0,1]$, $B(y,
c^{-j_n\rho})$ contains at least one of the $\widehat
I_{[j_n\rho],k}$'s. The conclusion follows.

\medskip

\noindent
(iii) If $\beta\le 1$, the fact that $\{(x_n,\lambda_n)\}_n$ forms almost
surely a weakly redundant system is a consequence of the estimates
obtained in the proofs of Lemma 5 and 8 of \cite{JAFFLEVY} for the
numbers $\widetilde N_{j,k}=\#\{n\in T_j:\ x_n\in
[k2^{-j},(k+1)2^{-j}]\}$.

If $\beta>1$, computations patterned after those performed in proving
(ii) show that if $\varepsilon\in (0,\beta-1)$, with probability 1,
there are infinitely many integers $j$ such that for all
$k\in\{0,\dots,c^{j}-1\}$, $\#\{n\in T_j:\ x_n\in I_{j,k}^c\}\ge
c^{j\varepsilon}$.
\end{proof}

\medskip 
{$\bullet$ Random family based on uniformly distributed points.}

\medskip
Let $\{x_n\}_n$ be a sequence of points independently and uniformly
distributed in $[0,1]^d$ and $\{\lambda_n\}_n$ a non-increasing
sequence of positive numbers.

We do not know conditions ensuring that
(\ref{P0r}) holds for some non-trivial measure $m$. The following
Proposition concerns (\ref{P0}) and weak redundancy.
\begin{proposition}
Let $\beta=\limsup_{j\to\infty}j^{-1}{\log_2 \# T_j}$.\\
{\em 1.}
Suppose that $\limsup_{n\ra+\infty} \left(\sum_{p=1}^n
{\lambda_p}/{2}\right) - d\, \log n = +\infty$. This implies $\beta\ge
1$. With probability 1 (\ref{recouvre1}) holds.

\smallskip
\noindent
{\em 2.}
Suppose that $\beta\le 1$. With probability 1, $\{(x_n,\lambda_n)\}_n$
is weakly redundant.
\end{proposition}
As a consequence, if $\lambda_n={\gamma}/{n}$ for some $\gamma>2d$
then, with probability 1, $\{(x_n,\lambda_n)\}_n$ is weakly redundant
and (\ref{recouvre1}) holds.
\begin{proof}
(i) It is Proposition 9 of \cite{KAHANE1}.

\smallskip

\noindent
(ii) The estimates of \cite{JAFFLEVY} invoked in the proof of
Proposition~\ref{bbb}(iii) also concern $\widehat N_{j,k}=\#\{n\in
T_j:\ x_n\in [k2^{-j},(k+1)2^{-j}]\}$ for the example we are dealing
with (i.e. $(x_n)$ is a sequence of i.i.d. uniform variables) when
$d=1$. In particular, when $d=1$, a sufficient condition for the
system to be weakly redundant is $ \beta\le 1$. Since a random
variable with uniform distribution in $[0,1]^d$ is a random vector in
$\mathbb{R}^d$ which components are independent uniform random
variables in $[0,1]$, the same property holds in dimension $d$ if
$\beta\le 1$.
\end{proof}

\subsection{Examples of measures $\mu$ and $m$, 
Interpretations of the property $\mathcal{P}^\rho_M$}
\label{sec_measures}

We give interpretations only for $\mathcal{P}^1_M$, since
$\mathcal{P}^\rho_M$ contains similar~information.

\smallskip

Given the measure $\mu$ and the exponent $\alpha>0$, there is
typically an uncountable family of values of $\beta>0$ such that
properties (\ref{P2}), (\ref{P1}), {\bf (3)} and {\bf (4)} of
Definition~\ref{deffubi} hold for many systems
$\{(x_n,\lambda_n)\}_n$. Consequently, one seeks for the largest value
of $\beta$. It follows from the study of the multifractal nature of
statistically self-similar (including the deterministic) measures we
deal with that, in general, this optimal value is given by
$\beta=\tau_\mu^*(\alpha)$ (see formulas (\ref{deftau}) and
(\ref{deftaue})).

\smallskip

We select four classes of measures to which Theorem \ref{theor_lower}
is applicable. Other examples can be found in
\cite{FAN2,BM1,BACRYMUZY,BM2,MOIBARRALREFINED}. We keep in mind item
{\em 3.} of Remark~\ref{rem}.

\smallskip

For the rest of this section the sequences $\{x_n\}_{n\in\mathbb{N}}$
and $\{\lambda_n\}_{n\in\mathbb{N}}$ are fixed, and we assume that
$(0,1)^d\subset \limsup_{n\to\infty}B(x_n,\lambda_n/2)$.

\smallskip

For $C,\kappa,r>0$ and $\gamma>1/2$, let $\varphi_C(r)=C
{|\log(r)|}^{-1/2}\big ({\log\log|\log (r)|}\big )^{1/2}$,
$\widetilde\varphi_\kappa(r)=\big (\log|\log(r)|\big )^{-\kappa}$, and
$\psi_{\gamma}(r)=C{|\log(r)|}^{-1/2}\big (\log |\log(r)|\big
)^{\gamma}$.


\medskip

{$\bullet$ Product of $d$ multinomial measures and
frequencies of digits}

\medskip

Let $(\pi^{(i)}_0,\dots,\pi^{(i)}_{c-1})$, $1\le i\le d$, be $d$
probability vectors with positive components such that
$\sum_{l=0}^{c-1}\pi^{(i)}_j=1$, $\forall \, 1\le i\le d$. For $1\le
i\le d$ let $\mu^{(i)}$ be the multinomial measure on $[0,1]$
associated with $(\pi^{(i)}_0,\dots,\pi^{(i)}_{c-1})$, and
$\mu=\mu^{(1)}\otimes\cdots\otimes\mu^{(d)}$ the product measure of
the $\mu^{(i)}$ on $[0,1]^d$. One has $ \tau_{\mu^{(i)}}(q)=-\log_c
\sum_{k=0}^{c-1}(\pi_k^{(i)})^q \ \mbox{ and } \ \tau_\mu
(q)=\sum_{i=1}^d\tau_{\mu^{(i)}}(q)$. It is convenient to take $\alpha
=\tau_\mu'(q)$ for some given $q\in \mathbb{R}$. Let us then define
$\beta=\tau_\mu^*(\alpha)=q\tau_\mu'(q)-\tau_\mu(q)$, and
$\mu_q=\mu_q^{(1)}\otimes\cdots\otimes\mu_q^{(d)}$, where
$\mu_q^{(i)}$ is the multinomial measure associated with the vector
$\big (c^{\tau_{\mu^{(i)}}(q)}(\pi_0^{(i)})^q,\dots,
c^{\tau_{\mu^{(i)}}(q)}(\pi_{c-1}^{(i)})^q\big)$.

It is proved in \cite{MOIBARRALREFINED0} that each measure $\mu^{(i)}$
satisfies properties (\ref{P2}), (\ref{P1}), {\bf (3)} and {\bf (4')}
with the exponents $\alpha_i=\tau'_{\mu^{(i)}}(q) $ and
$\beta_i=q\tau_{\mu^{(i)}}'(q)-\tau_{\mu^{(i)}}(q)$, and with $m$
equal to $\mu_q^{(i)}$.  This requires some work, because the masses
of the $c$-adic boxes and of their immediate neighbors need to be
controlled. One can choose $m^I\circ f_{I}^{-1}=m=\mu_q^{(i)}$, and
{\bf (3)} and {\bf (4')} do not matter. Moreover, $(\varphi,\psi)$ is
of the form $(\varphi_C,\psi_{\gamma})$.

Now, in terms of conditioned ubiquity, it is interesting to recall the
well-known interpretation of the conditions (\ref{P2}) and (\ref{P1}),
which hold for each $\mu^{(i)}$, in terms of $c$-adic expansions
(recall Section 1 and the definition (\ref{defphi}) of $\phi_{k,j}$):
For $\mu^{(i)}$-almost every point $x_i\in [0,1]$, for every $0\le
k\le c-1$, for all $y \in I_{j,k_{x_i}-1}\cup I_{j,k_{x_i}}\cup
I_{j,k_{x_i}+1}, \ \lim_{j\to\infty}
\phi_{k,j}(y)=c^{\tau_{\mu^{(i)}}(q)}(\pi_k^{(i)})^q$.

\smallskip

The previous remarks yield the following result, which implies
(\ref{JaBe}).

\begin{proposition}
Let $q\in\mathbb{R}$. The measure $\mu$ satisfies properties
(\ref{P2}), (\ref{P1}), {\bf (3)} and {\bf (4')} with $\alpha
=\tau_\mu'(q)$, $\beta=\tau_\mu^*(\alpha)$, $(\varphi,\psi)$ of the
form $(\varphi_C,\psi_{\gamma})$, and $m^I\circ f_I^{-1} =m=\mu_q$ for
all $ I\in \mathbf{I}$.

Moreover, there exists a sequence $\ep_n\searrow 0$ such that, when
applying Theorem~\ref{theor_lower}, property $
\mathcal{Q}(x_n,\lambda_n,1,\alpha,\varepsilon^1_{M,n})$ in
(\ref{defgene}) can be replaced by the following condition in terms of
$c$-adic expansion: for every $1\le i\le d$, for every $0\le k\le
c-1$, $\displaystyle \left|\phi_{k,[\log_c(\lambda_n^{-1})]}(x_{n,i})
-c^{\tau_{\mu^{(i)}}(q)}(\pi_k^{(i)})^q\right |\le \ep_n$, where
$x_n=(x_{n,1},\dots,x_{n,d})$.
\end{proposition}

{$\bullet$ Gibbs measures and
average of Birkhoff sums}

\medskip

Let $\phi$ be a $(1,\dots,1)$-periodic H\"older continuous function on
$\mathbb{R}^d$. Let $T$ be the transformation of $[0,1)^d$ defined by
$T\big ((x_1,\dots,x_d)\big )=(cx_1\mod 1,\dots,cx_d\mod 1)$. For
$k\in \mathbb{N}$, let $T^k$ denote the $k^{\mbox{th}}$ iteration of
$T$ ($T^0=\mbox{Id}_{[0,1)^d}$). For every $x\in [0,1)^d$ and $n\ge
1$, let us also define the $n^{\mbox{th}}$ Birkhoff sum of $x$,
$S_n(\phi)(x)=\sum_{k=0}^{n-1}\phi\big (T^k(x)\big )$ as well as
$D_n(\phi)(x)=\exp \big (S_n(\phi)(x)\big )$.

The Ruelle Perron-Frobenius theorem (see \cite{POLL}) ensures that the
probability measures $\mu_n$ given on $[0,1]^d$ by $ \mu_n
(dx)=D_n(\phi)(x)\, dx/{\int_{[0,1)^d}D_n(\phi)(u)\, du}$ converges
weakly to a probability measure $\mu$ which is a Gibbs state with
respect to the potential $\phi$ and the dynamical system
$([0,1)^d,T)$. The multifractal analysis of $\mu$ is performed in
\cite{FAN2,FAN5} for instance. With $\phi$ is also associated the
analytic function $\displaystyle L:\ q\in\mathbb{R}\mapsto d\log
(c)+\lim_{n\to\infty}{j}^{-1}\log \int_{[0,1)^d}D_n(q\phi)(u)\, du$,
which is the topological pressure of $q\phi$. One has $
\tau_\mu(q)=\frac{qL(1)-L(q)}{\log (c)}$. For $q\in\mathbb{R}$, let
$\mu_q$ be the Gibbs measure defined as $\mu$, but with the potential
$q\phi$.

Then, the structure of $\mu$ combined with the H\"older regularity of
$\phi$ and the law of the iterated logarithm (see Chapter 7 of
\cite{PHILIPPSTOUT}) yield

\begin{proposition}
Let $q\in\mathbb{R}$. The measure $\mu$ satisfies properties
(\ref{P2}), (\ref{P1}), {\bf (3)} and {\bf (4')} with $\alpha
=\tau_\mu'(q)$, $\beta=\tau_\mu^*(\alpha)$, both $\varphi$ and $\psi$
of the form $\varphi_C$, and $m^I\circ f_I^{-1} =m=\mu_q$ for all $
I\in \mathbf{I}$.

There exists $C>0$ such that, when applying Theorem~\ref{theor_lower},
in (\ref{defgene}) property $
\mathcal{Q}(x_n,\lambda_n,1,\alpha,\varepsilon^1_{M,n})$ can be
replaced in terms of average of Birkhoff sums by: $
\big|L'(q)-A_{[|\log_c(\lambda_{n})|]}(x_n) \big|\le
\varphi_C(\lambda_n)$, where $A_p(x):={S_p(\phi)(x)}/{p}$.
\end{proposition}


\medskip

{$\bullet$ Independent multiplicative cascades,
average of branching random~walks}

\medskip

For these random measures, the situation is subtle. Indeed, the study
achieved in \cite{MOIBARRALREFINED} concludes that property {\bf (4)}
can be satisfied for some systems $\{(x_n,\lambda_n)\}_{n\ge 1}$,
while the strong property {\bf (4')} fails because of the unavoidable
large values of $n_L$ for some $c$-adic boxes $L$.

Let us recall that these measures $\mu$ are constructed as
follows. Let $X$ be a real valued random variable. Let us define $L:\
q\in\mathbb{R}\mapsto d\log (c)+\log \mathbb{E}(e^{qX})$, and assume
that $L(1)<\infty$. For every $c$-adic box $J$ included in $[0,1]^d$,
let $X_J$ be a copy of $X$. Moreover, assume that the $X_J$'s are
mutually independent.  The branching random walk is then
\begin{equation}
\label{imc}
 \forall \, x\in [0,1)^d,\ \forall \, n\ge 1, \ S_n(x)=\sum_{J
    \in\mathbf{I} ,\ c^{-n}\le |J|\le c^{-1},\ x\in J} X_J.
\end{equation}

The measure $\mu$ is obtained as the almost sure weak limit of the
sequence $\mu_n$ on $[0,1]^d$ given by $\displaystyle \mu_n(dx)=\big
(\mathbb{E}(e^X)\big )^{-n}e^{S_n(x)}\, dx$.

Let $ \theta:\ q\in\mathbb{R}\mapsto \frac{qL(1)-L(q)}{\log (c)}$. In
\cite{MANDEL1,KP}, it is shown that $\theta'(1^-)>0$ is a necessary
and sufficient condition for $\mu$ to be almost surely a positive
measure with support equal to $[0,1]^d$. The multifractal nature of
$\mu$ or of variants of $\mu$ has been investigated in many works
\cite{K5,HoWa,FALC,O1,ArPa,Mol,BARRAL1}.  We need to consider the
interior $\mathcal{J}$ of the interval $\{q\in\mathbb{R}:\
\theta'(q)q-\theta(q)>0\}$.

\smallskip

For every $q\in \mathcal{J}$ and every $c$-adic box $I$ in $[0,1)^d$,
let us introduce the sequences of measures $\mu_{q,n}$ and $
m_{q,n}^I$ defined as follows: $\mu_{q,n}$ is defined as $\mu_n$ but
using $X_J(q):=qX_J$ instead of $X_J$ in (\ref{imc}), and $ m_{q,n}^I$
is defined as $\mu_{q,n}$ but with $qX_{f_I^{-1}(J)}$ instead of
$X_J(q)$ in (\ref{imc}).

It is shown in \cite{BARRAL1} that, with probability 1, $\forall\,
q\in \mathcal{J}$, the measures $\mu_{q,n}$ converge weakly to a
positive measure $\mu_q$ on $[0,1]^d$; In addition, $\forall \, q\in
\mathcal{J}$, for every $c$-adic box $I$ of generation $\ge 1$, the
sequence of measures $ m_{q,n}^J$ converges weakly to a measure $
m_{q}^I$ on $[0,1]^d$, and $\tau_\mu(q)=\theta(q)$ on~$\mathcal{J}$.

The following result is a consequence of Theorem 4.1 in
\cite{MOIBARRALREFINED}.

\begin{proposition}\label{propmand}
Suppose that $\limsup_{n\to\infty}B(x_n,\lambda_n/4)\supset (0,1)^d$.

For every $q\in \mathcal{J}$, with probability 1 (and also with
probability 1, for almost every~$q\in \mathcal{J}$), $\mu$ satisfies
properties (\ref{P2}), (\ref{P1}), {\bf (3)} and {\bf (4)} with the
exponents $\alpha=\tau'_{\mu}(q) $ and $\beta=\tau_\mu^*(\alpha)$,
$(\varphi,\psi)$ of the form $(\widetilde\varphi_\kappa,\psi_\gamma)$,
$m=\mu_q$, $m^I\circ f^{-1}_I= m_{q}^I$ for all $ I\in \mathbf{I}$,
and $\mathcal{D}=\mathbb{Q}\cap (1,\infty)$.

There exists $\gamma>1/2$ such that, when applying
Theorem~\ref{theor_lower}, in (\ref{defgene}) property $
\mathcal{Q}(x_n,\lambda_n,1,\alpha,\varepsilon^1_{M,n})$ can be
replaced in terms of average of branching random walks by: $
\big|L'(q)-A_{[|\log_c(\lambda_{n})|]}(x_n)
\big|\le\psi_\gamma(2\lambda_n)$, where $A_p(x):={S_p(x)}/{p}$.
\end{proposition}

\medskip

{$\bullet$ Poisson cascades and average of
covering numbers in the case $d=1$.}

\medskip
 Let $\xi>0$ and $S$ a Poisson
point process in $\mathbb{R}\times (0,1)$ with intensity $\Lambda$
given~by $\displaystyle \Lambda (ds\, d\lambda)={\xi}{ds
d\lambda}/2{\lambda^2}$. For every $c$-adic box $I$ of $[0,1]$, define
$S_I=\big \{(f_I^{-1}(t), |I|^{-1}\lambda): (t,\lambda)\in S,\ \lambda
<|I|\big\}$. The point process $S_I$ is a copy of $S$.

\smallskip

For every $t\in [0,1]$ and $\ep\in (0,1]$, the covering number of $t$
at height $\ep$ by the Poisson intervals $\{(s-\lambda,s+\lambda):
(s,\lambda)\in S\}$ is defined by
\begin{eqnarray*}
N^S_\ep (t)=\!\!\!\!\!\sum_{(t,\lambda)\in S, \ \lambda\ge
\ep}\!\!\!\!\!\mathbf{1}_{\{(s-\lambda,s+\lambda)\}}(t)=\#\big
\{(s,\lambda)\in S:\ \lambda\ge \ep, \ t\in (s-\lambda,s+\lambda)\big
\}.
\end{eqnarray*}
The measure $\mu$ on $[0,1]$ is the almost sure weak limit, as $\ep\to
0$, of
\begin{equation}
\label{cpc}
\mu_\ep(dt)= \big (\mathbb{E}(e^{N^S_\ep (t)}\big )\big )^{-1}e^{N^S_\ep (t)}\,
dt =  \ep^{\xi (e-1)}e^{N^S_\ep (t)}\, dt.
\end{equation}
Let $L:q\in\mathbb{R}\mapsto {\xi}^{-1}+e^q-1$, and let $\theta:\
q\in\mathbb{R}\mapsto \xi \big (qL(1)-L(q)\big )$.

In \cite{BM1}, it is shown that $\theta'(1^-)>0$ is a necessary and
sufficient condition for $\mu$ to be almost surely a positive measure
supported by $[0,1)^d$.  Let $\mathcal{J}=\{q\in\mathbb{R}:\
\theta'(q)q-\theta(q)>0$.  It is also shown in \cite{BM1} that, with
probability 1, for all $q\in \mathcal{J}$, the measures $\mu_{q,\ep}$
on $[0,1]$ given by $\mu_{q,\ep}(dt)=\ep^{\xi(e^q-1)}e^{qN^S_\ep
(t)}\, dt$ converge weakly, as $\ep\to0$, to a positive measure
$\mu_q$ on $[0,1]$; moreover, for every $q\in \mathcal{J}$, for every
$c$-adic interval $I$ of generation $\ge 1$, the family of measures $
m_{q,\ep}^I$ constructed as $\mu_{q,\ep}$ but with $N^{S_I}_\ep (t)$
instead of $N^{S}_\ep (t)$ in (\ref{cpc}) converges weakly, as $\ep\to
0$, to a measure $m_{q}^I$ on $[0,1]$; finally, one has
$\tau_\mu(q)=\theta(q)$ on $\mathcal{J}$.

The same conclusions as in Proposition~\ref{propmand} hold if $
\mathcal{Q}(x_n,\lambda_n,1,\alpha,\varepsilon^1_{M,n})$ is replaced
by $\Big|L'(q)+ \frac{1}{\xi\log (\lambda_n)} N_{\lambda_n}(x_n)\Big|
\le \psi_\gamma(2\lambda_n)$.

More on covering numbers and related questions can be found in
\cite{BARRALFAN,BARRALFAN2}.

\subsection{Example where $\dim \,\big (
\limsup_{n\to\infty}B(x_n,{\lambda_n}/{2})\big )<d$.}
\label{sec_related}

Let us return to the example of Gibbs measures $\mu$ in
Section~\ref{sec_measures}. Let $q_0>0$. Fix $\mathcal K$ a subset of
$\mathbb{R}$ such that $\tau_\mu'(\mathcal K)\cap (\tau_\mu'(q_0),
\tau_\mu'(-q_0))=\emptyset$. Define the~system
\[\{(x_n,\lambda_n)\}= \Big\{\left(\left
(\mathbf{k}+{\mathbf{1}}/{\mathbf{2}}\right )c^{-j},
c^{-j}\right):\displaystyle \frac{\log \mu\left(B\left(\left
(\mathbf{k}+{\mathbf{1}}/{\mathbf{2}}\right ),
c^{-j}\right)\right)}{-j\log (c)}\in \mathcal K \Big \}.
\]
Let $S=\limsup_{n\to\infty}B(x_n,{\lambda_n}/{2})$. For every $q\in
\mathcal K$, one has $\mu_q (S)=1$ and $\dim\, S\le \max \big
(\tau_\mu^*(\tau_\mu'(-q_0)), \tau_\mu^*(\tau_\mu'(q_0))\big )<d$.

\bibliographystyle{plain}


%



%
%

\end{document}